\documentclass[11pt,a4paper,rotate,colorlinks]{article}
\RequirePackage[OT1]{fontenc}
\RequirePackage{amsthm,amsmath}
\RequirePackage[colorlinks,citecolor=blue,urlcolor=blue]{hyperref}
\RequirePackage{natbib}
\RequirePackage{hypernat}

\usepackage{amsfonts}
\usepackage{rotating}
\usepackage{hyperref}
\usepackage[english]{babel}
\usepackage{latexsym}
\usepackage{amsmath}
\usepackage{amssymb,bbm}
\usepackage{supertabular}
\usepackage{setspace}
\usepackage{xspace}
\usepackage{url}
\usepackage{color}
\usepackage{graphicx}
\usepackage{listings}
\usepackage{natbib}

\newtheorem{assumption}{Assumption}
\newtheorem{lemma}{Lemma}[section]
\newtheorem{theorem}{Theorem}[section]
\newtheorem{corollary}{Corollary}[section]
\newtheorem{remark}{Remark}

 \normalfont
 \makeatletter

\def \sgn {\mbox{sign}}
\def\ar{\textsc{AR }}

\def\dar{\textsc{DAR }}

\def \st {\mbox{st}}

\begin{document}
\title{Strict Stationarity Testing and GLAD Estimation  of Double Autoregressive Models}
\author{Shaojun Guo$^{a}$,\,\, Dong Li$^{b}$, \,\, Muyi Li$^{c,d}\thanks{%
Corresponding author: Muyi Li, Economics Building, Xiamen
University, Xiamen 361005, Fujian Province, China. Email: limuyi@xmu.edu.cn. Tel.:
+86-592-2182387. Fax: +86-592-2187708.}$ \\
$^{a}$ Institute of Statistics and Big Data, Renmin University of China, Beijing 100872, China.\\
$^{b}$ Center for Statistical Science and Department of Industrial Engineering, \\ Tsinghua University, Beijing 100084, China. (malidong@tsinghua.edu.cn)\\
$^{c}$ MOE Key Laboratory of Econometrics,\\ Wang Yanan Institute for Studies in Economics (WISE), Xiamen University,\\
$^{d}$ Department of Statistics, School of Economics,
Xiamen University.\\ (limuyi@xmu.edu.cn)}
\date{}
\maketitle

\begin{abstract}
In this article we develop a tractable procedure for testing strict stationarity in a double autoregressive  model and formulate the problem as testing if the top Lyapunov exponent is negative.  Without strict stationarity assumption, we construct a consistent estimator of the associated top Lyapunov exponent and employ a random weighting approach for its variance estimation, which in turn are used in a $t$-type test.  We also propose a GLAD estimation for parameters of interest, relaxing key assumptions on the commonly used QMLE. All estimators, except for the intercept, are shown to be consistent and asymptotically normal in both stationary and explosive situations. The finite-sample performance of the proposed procedures is evaluated via Monte Carlo simulation studies and a real dataset of interest rates is analyzed.

\medskip


\end{abstract}

\section{Introduction}
The assumption of strict stationarity is pivotal in nonlinear time series inference and forecasting. Testing stationarity in the context of linear time series models has been well documented, such as various unit root tests. 
However, this testing problem may pose considerably more challenges in a nonlinear setting. Recently, \cite{fz12,fz13} considered strict stationarity testing for GARCH models (\cite{engle}, \cite{bollerslev}) based on the sign of the associated top Lyapunov exponent. To our best knowledge, they are the first to address such a testing issue under the GARCH framework.  
 
Besides GARCH, double autoregressive (DAR) model  is another important conditional heteroscedastic one.  The first-order \dar (hereafter \dar(1)) model is defined as
\begin{eqnarray}\label{dar}
y_t=\phi_0 y_{t-1}+ \eta_t\sqrt{\omega_0+\alpha_0 y_{t-1}^2}, \quad t=0, \pm1, \pm2,...,
\end{eqnarray}
where $\phi_0\in R$, $\omega_0>0$, $\alpha_0 > 0$, $\{\eta_t\}$ is a sequence of independent and identically distributed (i.i.d.) random variables and independent of $\{y_{j}, j<t\}$.
Model (\ref{dar}) is a special case of
the \textsc{ARMA-ARCH} models in \cite{weiss} and of the nonlinear \ar models in \cite{cpu}. It is different from Engel's \textsc{ARCH} model when $\phi_0\neq 0$. There has been considerable work  on \dar models and its generalizations, see, e.g., \cite{tsay}, \cite{lu}, \cite{ling04, ling07}, \cite{zhu}, \cite{NR2014}, \cite{li15a},
\cite{li15b}, \cite{liguodong},  \cite{Zhu:Zheng:Li}.

\begin{figure}[!htbp]
\begin{center}
\includegraphics[width = 70mm, height = 80mm]{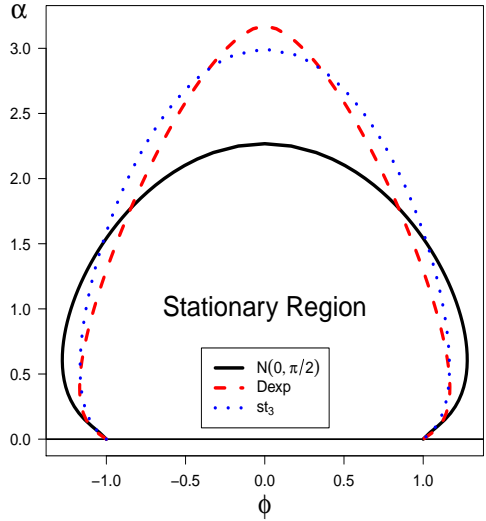}
\caption{The (strictly) stationary region determined by $\{(\phi, \alpha): E\log|\phi+\eta\sqrt{\alpha}|<0\}$ for
$\eta\sim N(0, \pi/2)$, the double exponential distribution with density $f(x)=0.5\exp(-|x|)$, and the standardized
Student's $t_3$ $(\mathrm{st}_3)$ with density $f(x)=4\pi^2/(\pi^2+4x^2)^2$.}
\label{boundary}
\end{center}
\end{figure}

The strict stationarity condition for model (\ref{dar}) has been well formulated and is closely related to the top Lyapunov exponent
\begin{flalign*}
\gamma_0 = E\log|\phi_0+\eta\sqrt{\alpha_0}|,
\end{flalign*}
where $\eta$ is a generic random variable with the same distribution as $\eta_t$. \cite{bk} proved that $\gamma_0<0$ is sufficient for strict stationarity of model (\ref{dar}),
while \cite{cll} proved that it is (almost) necessary. Thus, testing strict stationarity is equivalent to testing the following hypothesis:
\begin{flalign}\label{testone}
H_0:~\gamma_0<0\quad \mathrm{v.s.} \quad H_1:~\gamma_0  \ge 0.
\end{flalign}
To test (\ref{testone}), a challenge is  to obtain a good estimator of $\gamma_0$ under the null and alternative hypotheses. Surprisingly, the majority of literature is focused on the inference of $(\phi_0, \alpha_0, \omega_0)$ rather than $\gamma_0$. \cite{cp2005} discussed the estimation of $\gamma_0$  in the stationary case via simulation studies so as to check the stationarity condition.  Nevertheless, they do not consider its asymptotic properties. This may be partly due to the fact that the estimation of $\gamma_0$ is nonstandard even in the stationary situation and the related asymptotic theory is hard to pursue. Fig. \ref{boundary} shows the strictly stationary region of model (\ref{dar}) for three different distributions of $\eta$, indicating that $\gamma_0$ depends not only on the intrinsic parameters but also on the distribution of the underlying innovations.

It is notable that \cite{ling04} considered the following testing problem for model (\ref{dar}):
\begin{flalign}\label{test2}
H_0: ~(\phi_0,\alpha_0) = (\pm 1,0) \quad \mathrm{v.s.} \quad H_1:~(\phi_0,\alpha_0) \neq (\pm 1,0).
\end{flalign}
Under  $H_0$ in (\ref{test2}), $\{y_t\}$ is a standard unit root process. However, we cannot claim that $\{y_t\}$ is stationary even if $H_0$ in (\ref{test2}) is rejected. \cite{ling04} further points out that testing $\gamma_0 <0$ in (\ref{testone}) remains an interesting issue. One aim of this paper is to develop a tractable procedure for testing (\ref{testone}) within the DAR framework.

This paper has two major contributions. First, we propose a data-driven procedure for testing (\ref{testone}). The basic steps are as follows. We first provide a natural estimator of $\gamma_0$ and demonstrate that it possesses nice asymptotic properties without any stationarity assumption. We then develop a random weighting method to approximate its asymptotic covariance. A feature of our procedure is that it does not rely on the  expression of asymptotic variance and turns out to capture the sampling uncertainty adaptively. Based on these estimators, we propose a $t$-type statistic for testing strict stationarity and construct a consistent critical region. It is remarkable that our basic idea is similar to that in \cite{fz12} but the underlying large-sample theories own considerably more mathematical gaps. In other words, the proof techniques in \cite{fz12} fail to apply here. To solve this problem, we rely on modern empirical process theory. A new theoretical insight on this issue is provided in Section \ref{threesect}.

In the literature, numerous methods have been developed for inference on  model (\ref{dar}).  \cite{ling04,ling07} considered the quasi-maximum likelihood estimation (QMLE) of $(\phi_0,\alpha_0,\omega_0)$ and proved its asymptotics under the  conditions $\gamma_0<0$ and $E\eta^4 < \infty$. On the other hand, in the explosive case ( i.e., $\gamma_0 > 0$), \cite{lili} 
PDF Rendering Error Something went wrong while rendering this PDF.
 investigated  a constrained QMLE of $(\phi_0, \alpha_0)$ in the sense that $\omega_0$  is fixed and obtained asymptotic normality when $\eta$ is standard normal.
Recently, \cite{cll} studied an unconstrained QMLE of $(\phi_0, \alpha_0)$ when $\eta$ is symmetric with $E\eta^4<\infty$.
One disadvantage of the QMLE is that the assumption of $E\eta^4 < \infty$ is indispensable for valid inference but too restrictive in practice. To tackle this challenge, a robust estimation procedure is provided, see \cite{cp2005} and \cite{zhu} for the stationary case.
In the presence of nonstationarity or heavy-tailed noises, the inference becomes more challenging and
no results are available so far. 

The second contribution of this article is to offer a unified framework for parameter estimation of model (\ref{dar}) in both stationary and explosive cases. This framework does not require the fourth moment of the innovation $\eta_t$ to be finite and hence can cover heavy-tailed cases. Specifically,
we propose an unconstrained global least absolute deviation estimation (GLADE) for $(\phi_0,\alpha_0,\omega_0)$ when $\omega_0$ is not fixed. Here, `global' means that the convergence rate of the estimator is first obtained and the limiting distribution is then derived, see \cite{zhling,zhu}.
Under mild conditions, the proposed estimator of $(\phi_0,\alpha_0)$ are always strongly consistent and asymptotically normal when the data mechanism is stationary or explosive. To measure the accuracy of the estimator, we again propose a random weighting method to estimate its asymptotic covariance matrix. It is worth stressing that even if the GLADE of $(\phi_0,\alpha_0)$ is consistent in every situation, the intercept term $\omega_0$ is only consistent in the stationary case. This further demonstrates the importance of testing the sign of $\gamma_0$.

The remainder of the paper is organized as follows. Section \ref{twosect} presents the global LADE with asymptotic properties and discusses how to estimate the asymptotic variance matrix via the random weighting approach. In Section \ref{threesect}, we address the estimation of $\gamma_0$ and tests of stationarity and nonstationarity. The relevant asymptotics are also established.
Section \ref{foursect} reports numerical results on the finite-sample performance of the proposed methodology and analyzes a real dataset of interest rates. A concluding remark is in Section \ref{fivesect} and all the technical proofs can be found in the Appendix.

\section{Global Least Absolute Deviation Estimation}\label{twosect}
\subsection{Asymptotic Properties}
Suppose that the observations $\{y_0, y_1,...,y_n\}$ are from model (\ref{dar}). When $\eta$ is double exponential, the log-likelihood
function (ignoring a constant) can be written as
\begin{flalign*}
L_n(\theta)=\sum_{t=1}^n\bigg\{\frac{1}{2}\log(\omega+\alpha y_{t-1}^2)+\frac{|y_t-\phi y_{t-1}|}{\sqrt{\omega+\alpha
y_{t-1}^2}}\bigg\},
\end{flalign*}
where $\theta=(\phi, \alpha, \omega)^T$ is the parameter. The proposed estimator is defined as
\begin{eqnarray}\label{lad}
\hat{\theta}_n=(\hat{\phi}_n, \hat{\alpha}_n, \hat{\omega}_n)^T=\arg\min_{\theta\in\Theta}L_n(\theta),
\end{eqnarray}
where $\Theta$ is a compact subset of $R\times R_{+}^2$ containing the true value
$\theta_0=(\phi_0, \alpha_0, \omega_0)^T$. Here $R_+=(0, \infty)$. Since we do not assume that $\eta$ is double exponential, the estimator $\hat{\theta}_n$ is often called the quasi-maximum exponential likelihood estimator as in \cite{zhu} or  the least absolute deviation estimator (LADE) as in \cite{cp2005}. Throughout the paper,
following a traditional naming convention in literature,
we refer to $\hat{\theta}_n$ as the (global) LADE of $\theta_0$.

\begin{theorem}\label{consistency}
Suppose that $\{\eta_t\}$ is i.i.d. and symmetric  with $E|\eta_t|=1$ and $\Theta$ is compact.
Then, for the DAR(1) model (\ref{dar}), the LADE defined in (\ref{lad}) satisfies the following properties.\\
$(\mathrm{i}).$ If $\gamma_0<0$, then $\hat{\phi}_n\rightarrow \phi_0$, $\hat{\alpha}_n\rightarrow\alpha_0$, and $\hat{\omega}_n\rightarrow\omega_0$ a.s. as $n\rightarrow\infty$.\\
$(\mathrm{ii}).$ If $\gamma_0>0$, then $\hat{\phi}_n\rightarrow \phi_0$ and $\hat{\alpha}_n\rightarrow\alpha_0$ a.s.
as $n\rightarrow\infty$.
\end{theorem}

To further discuss asymptotic distribution of $\hat{\theta}_n$,  we need three assumptions.
\begin{assumption}\label{asm1}
$\{\eta_t\}$'s are $i.i.d.$ with $E|\eta_1|=1$ and $\kappa_\eta= E\eta_1^2<\infty$.
\end{assumption}
\begin{assumption}\label{asm2}
The density $f(x)$ of $\eta_1$ is symmetric and bounded continuous on $R$ with $f(0)>0$.
\end{assumption}
\begin{assumption}\label{asm3}
$\Theta$ is compact and $\theta_0$ is an interior point of $\Theta$.
\end{assumption}

\begin{theorem}\label{theorem2}
Suppose that Assumptions \ref{asm1}--\ref{asm3} hold.\\
$(\mathrm{i})$. If $\gamma_0<0$, then
\begin{flalign*}
\sqrt n(\hat{\theta}_n-\theta_0)=O_p(1)\quad\mbox{and}\quad\sqrt n(\hat{\theta}_n-\theta_0)\stackrel{\mathcal{L}}{\longrightarrow}N\big(0,~\mathcal{J}_S\big)
\end{flalign*}
as $n \to \infty$, where `$\stackrel{\mathcal{L}}{\longrightarrow}$' stands for the convergence in distribution,
\begin{flalign}\label{w1}
\mathcal{J}_S=\mathrm{diag}\left\{\frac{1}{4\sigma_{11}f^2(0)},~4(\kappa_\eta-1) \Sigma^{-1}\right\},
\end{flalign}
and both $\Sigma$ and $\sigma_{11}$ are defined as
\begin{flalign}\label{sigma}
\Sigma=\left(\begin{array}{cc}
\sigma_{22} &\sigma_{12}\\
\sigma_{12} &\sigma_{02}
\end{array}
\right)\quad \mbox{with} \quad \sigma_{ij}=E\left\{\frac{y_t^{2i}}{(\omega_0+\alpha_0y_t^2)^j}\right\},
\quad i,j = 0,1,2.
\end{flalign}
$(\mathrm{ii})$. If $\gamma_0>0$, then
\begin{flalign*}
\sqrt n(\hat{\phi}_n-\phi_0, \hat{\alpha}_n-\alpha_0)^T=O_p(1)\quad \mbox{and}\quad
\sqrt n(\hat{\phi}_n-\phi_0, \hat{\alpha}_n-\alpha_0)^T\stackrel{\mathcal{L}}{\longrightarrow}
N\Big(0, \mathcal{J}_N\Big),
\end{flalign*}
as $n \to \infty $, where
\begin{flalign}\label{w2}
\mathcal{J}_N = \mathrm{diag}\Big\{\frac{\alpha_0}{4f^2(0)},~4(\kappa_\eta-1)\alpha^2_0\Big\}.
\end{flalign}
\end{theorem}

\begin{remark}
From Theorem \ref{theorem2}, we can see that our LADE is global, i.e., the convergence rate is first
obtained and the limiting distribution then derived, see \cite{zhling, zhu}. This is totally different from the local LADE of regression or time series models.
On the other hand,
it is worth noting that $\hat{\omega}_n$ is inconsistent when $\gamma_0>0$. A similar
phenomenon was observed by \cite{cll}, who studied an unconstrained QMLE of an explosive DAR(1) model,
and by \cite{fz12} who
studied the QMLE of nonstationary GARCH(1,1) models, see also \cite{jra,jrb}.
\end{remark}

\begin{remark}
We explore  the hidden relationship between expressions (\ref{w1}) and (\ref{w2}),
which reveals why several results in the stationary situation can still be applicable in the explosive one.
Note that all $\sigma_{ij}$'s in (\ref{sigma}) are finite constants when $\gamma_0<0$ since $\{y_t\}$ is strictly stationary
and ergodic. In the explosive situation, i.e. $\gamma_0 > 0$, however,
$|y_t|$ diverges to infinity at an exponential rate as $t \rightarrow \infty$, see Theorem 1 in \cite{cll}
or Theorem 2.1 in \cite{liuf}.
Thus, $\sigma_{11}\rightarrow1/\alpha_0$, $\sigma_{22}\rightarrow1/\alpha_0^2$  and other $\sigma_{ij}$'s go to zero a.s. as $t\rightarrow\infty$, which implies that $\mathcal{J}_S$ can be reduced to
$\mathrm{diag}(\mathcal{J}_N, 0)$.
\end{remark}

\begin{remark}
Theorem \ref{theorem2} rules out the boundary case $\alpha_0=0$ or $\omega_0=0$.  Clearly, if
$\alpha_0=0$, the GLADE $\hat{\alpha}_n$ is not asymptotically normal at all since $\hat{\alpha}_n\geq0$.
A similar phenomenon can be found in  Example 5 for an ARCH(1) model in \cite{fz07}.
When $\omega_0=0$, model (\ref{dar}) reduces to a double AR(1) model without intercept, 
which has been studied by \cite{lgz2017}.
\end{remark}

\begin{remark}
Unlike the GARCH(1, 1) model in \cite{fz12}, we don't establish  asymptotic properties of the GLADE for the case  $\gamma_0=0$ since the asymptotic behavior of $|y_t|$ is unknown in this scenario. Although we may pay some extra cost, such as more unverifiably restrictive conditions on $\eta_t$ and $y_t$,  to obtain asymptotics of the GLADE, we don't pursue this at this moment and left it for future work.
\end{remark}

\begin{remark}
The GLADE proposed here can be easily extended to a strictly stationary DAR($p$) ($p\geq 1$) model with an intercept in the conditional mean. For a nonstationary DAR(1) model with an intercept in the conditional mean, it can be shown that the GLADE of this intercept is not consistent, just like $\omega_0$ in this case. However, for nonstationay DAR($p$) model with $p\geq 2$, we cannot establish its asymptotics following the routine of nonstationary high-order GARCH models as they pose considerably huge gap between these two different processes. We left this for future work.
\end{remark}

\begin{remark}
In standard linear time series analysis, nonstationarity
leads to completely different asymptotic inferences, such as the Dickey-Fuller statistics. However, this is not so in DAR models, where a unit root in the conditional mean does not necessarily imply nonstationarity as lagged observations of the process enter the conditional variance.
This point can be seen from Fig.\ref{boundary}. See also \cite{bk}, \cite{ling04} and \cite{NR2014}.
\end{remark}

\subsection{Variance Estimation}
Theorem \ref{theorem2} demonstrates that the proposed global robust estimator shares nice asymptotic normality in both cases. To evaluate the accuracy of the estimator, we need a consistent estimator of their asymptotic covariance matrix. Unlike the QMLE, however, their asymptotic covariance
matrix involves the density function $f(\cdot)$ and cannot be precisely estimated using the plug-in rules.
One choice is to use the nonparametric kernel method to estimate the density, see \cite{zhu}.
However, this method involves the choice of the bandwidth, which is another problem. To avoid density estimation, we propose a random weighting approach for covariance matrix estimation. This approach, as a variant of
the traditional wild bootstrap in \cite{wu1986}, was initially proposed by \cite{Jin2001}. So far, it has been widely used for statistical inference in regression and time series models when the asymptotic covariance matrix of certain estimator can not be directly estimated; see \cite{koul1994}, \cite{koul2002}, \cite{chen2008, chen2010}, \cite{lgd2014} and \cite{Zhu2016}, etc.

To be specific, let $\varpi _1,\cdots,\varpi_n$ be a sequence of
i.i.d. nonnegative random variables, with both mean and variance
equal to one. For example, the standard exponential distribution
satisfies this requirement. Define
\begin{flalign}\label{ladd}
L_n^*(\theta)=\sum_{t=1}^n \varpi _t\bigg\{\frac{1}{2}\log(\omega+\alpha y_{t-1}^2)+\frac{|y_t-\phi y_{t-1}|}{\sqrt{\omega+\alpha
y_{t-1}^2}}\bigg\}
\end{flalign}
and $\hat{\theta}_n^* = \arg \min_{\theta \in \Theta} L_n^*(\theta)$.
In the stationary case, the distribution of $\sqrt{n}(\hat{\theta}_n -\theta_0)$ can be approximated by the resampling distribution of $\sqrt{n}(\hat{\theta}_n^* - \hat{\theta}_n)$. It turns out that, in the explosive case, the distribution of $\sqrt{n}(\hat{\phi}_n - \phi_0, \hat{\alpha}_n - \alpha_0)^T$ can still be approximated by that of $\sqrt{n}(\hat{\phi}_n^* - \hat{\phi}_n, \hat{\alpha}_n^* - \hat{\alpha}_n)^T$.

\begin{theorem}\label{theorem3}
Suppose Assumptions \ref{asm1}--\ref{asm3} hold. \\
$(\mathrm{i})$. If $\gamma_0 < 0$, then, conditionally on the data $\{y_0,y_1,\cdots, y_n\}$, in probability,
\begin{flalign*}
\sqrt n(\hat{\theta}_n^* - \hat{\theta}_n) \stackrel{\mathcal{L}}{\longrightarrow} N\big(0,~\mathcal{J}_S\big),\quad \mbox{as $n \to \infty$;}
\end{flalign*}
$(\mathrm{ii})$. If $\gamma_0>0$, then, conditionally on the data $\{y_0,y_1,\cdots, y_n\}$, in probability,
\begin{flalign*}
\sqrt n(\hat{\phi}_n^* - \hat{\phi}_n, ~\hat{\alpha}_n^*- \hat{\alpha}_n)^T\stackrel{\mathcal{L}}{\longrightarrow}
~N\big(0, ~\mathcal{J}_N\big),\quad \mbox{as $n \to \infty$.}
\end{flalign*}
\end{theorem}

The proof of Theorem \ref{theorem3} is provided in the Appendix. The inference
procedure via random weighting is done as follows. First, i.i.d. $\{\varpi_1,\cdots,\varpi_n\}$ are
generated $B$ times from standard exponential distribution, where $B$ is a large number.
Each time, the minimizer $\hat{\theta}^*_n$ is computed.
Denote them as $\hat{\theta}^{*1}_n,\cdots,\hat{\theta}^{*B}_n$.
Then, the asymptotic variances of $\sqrt{n}(\hat{\phi}_n - \phi_0)$ and
$\sqrt{n}(\hat{\alpha}_n - \alpha_0)$ can be approximated by, respectively,
\begin{eqnarray*}
\hat{\sigma}^2_{\phi} = \frac{1}{B-1}\sum_{b = 1}^B \big(\hat{\phi}^{*b}_n - \bar{\phi}_n^* \big)^2 \ \ \mbox{and} \ \
 \hat{\sigma}^2_{\alpha} = \frac{1}{B-1}\sum_{b = 1}^B \big(\hat{\alpha}^{*b}_n - \bar{\alpha}_n^* \big)^2,
\end{eqnarray*}
with $\bar{\phi}_n^* = B^{-1}\sum_{b=1}^B \hat{\phi}^{*b}_n$ and $\bar{\alpha}_n^* = B^{-1}\sum_{b=1}^B \hat{\alpha}^{*b}_n$.
In the stationary case, the asymptotic variance of $\sqrt{n}(\hat{\omega}_n - \omega_0)$ is estimated analogously.

\section{Strict Stationarity Testing}\label{threesect}
As described in the Introduction, the parameter $\gamma_0$ plays a key role in characterizing the strict stationarity of DAR(1) model. It also determines the consistency of $\hat{\omega}_n$. In this section, we propose a consistent estimator $\hat{\gamma}_n$ of $\gamma_0$ and then construct $t$-type tests based on $\hat{\gamma}_n$ for whether $\{y_t\}_{t=1}^n$ is strictly stationary or not. The asymptotic normality of $\hat{\gamma}_n$ is established and  a random weighting approach is correspondingly proposed to estimate the asymptotic variance of $\hat{\gamma}_n$ in Section \ref{subsectone}. Section \ref{subsecttwo} shows that the $t$-type stationarity tests are consistent.

\subsection{Consistent Estimator of $\gamma_0$}\label{subsectone}
Define the rescaled residuals as
 \begin{flalign}\label{residual}
\hat{\eta}_t= \eta_t(\hat{\theta}_n), \quad \eta_t(\theta)=\frac{y_t-\phi y_{t-1}}{\sqrt{\omega+\alpha y_{t-1}^2}}.
 \end{flalign}
By the definition of $\gamma_0$ and the symmetry of $\eta_t$, a natural estimator of $\gamma_0$ is
\begin{flalign}\label{naturalgamma}
\tilde{\gamma}_{n} =\frac{1}{2n} \sum_{t=1}^n\big(\log|\hat{\phi}_n+\hat{\eta}_t\sqrt{\hat{\alpha}_n}|  + \log|\hat{\phi}_n - \hat{\eta}_t\sqrt{\hat{\alpha}_n}| \big).
\end{flalign}
To facilitate our proofs, we propose a truncated estimator for $\gamma_0$:
\begin{flalign}\label{truncategamma}
\hat{\gamma}_{n}=\frac{1}{2n} \Big (\sum_{t \in \mathcal{A}_1}\log|\hat{\phi}_n+\hat{\eta}_t\sqrt{\hat{\alpha}_n}| + \sum_{t \in \mathcal{A}_2}\log|\hat{\phi}_n -\hat{\eta}_t\sqrt{\hat{\alpha}_n}|\Big ),
\end{flalign}
where  $\mathcal{A}_1 = \{ t: \hat{\phi}_n + \hat{\eta}_t\sqrt{\hat{\alpha}_n} \in \mathcal{I}_n, 1 \le t \le n \} $ and $ \mathcal{A}_2 = \{t: \hat{\phi}_n - \hat{\eta}_t\sqrt{\hat{\alpha}_n} \in \mathcal{I}_n, 1 \le t \le n \}$
with $\mathcal{I}_n = [- n^2,-n^{-2}] \cup [n^{-2},n^2]$. Indeed, in practice, both estimates are almost identical unless $|\hat{\phi}_n + \hat{\eta}_t\sqrt{\hat{\alpha}_n}|$ or $|\hat{\phi}_n - \hat{\eta}_t\sqrt{\hat{\alpha}_n}|$ is extremely small for some $t$.

To prove the asymptotics of $\hat{\gamma}_n$, we impose the following assumption.
\begin{assumption} \label{asm4}
The density $f(\cdot)$ of $\eta_1$ is positive and differentiable a.s. on $R$ with
$\sup_{x \in R}f(x) < \infty $, and
$\int (\log|\phi_0 + x\sqrt{\alpha_0}|)\{f(x) + |xf'(x)|\}dx$
exists and  is finite.
\end{assumption}

\begin{remark}
The most commonly used distributions, like normal and Student's $t$ distributions, satisfy this condition. The double exponential (i.e., Laplace) distribution is also allowed since  it is differentiable except at one point.
\end{remark}

To state the asymptotic normality of $\hat{\gamma}_n$, we first introduce some notation. Denote
\begin{flalign*}
\mu_1 &=- \alpha_0^{-1/2}\int \log|\phi_0 + x\sqrt{\alpha_0}|f'(x)dx, \\
\mu_2 &=- \alpha_0^{-1/2}\int \log|\phi_0 + x\sqrt{\alpha_0}| \{f(x) + xf'(x)\} dx.
\end{flalign*}
For each $t$, define
\begin{flalign*}
\zeta_t &=(\zeta_{1t} + \zeta_{2t})/2, \\
\zeta_{1t} &= \log|\phi_0+ \eta_t \sqrt{\alpha_0}| - E (\log|\phi_0+\eta_t \sqrt{\alpha_0}|), \\
\zeta_{2t} &= \log |\phi_0-\eta_t \sqrt{\alpha_0}| - E (\log |\phi_0-\eta_t \sqrt{\alpha_0}|).
\end{flalign*}
For $\gamma_0 <0$, denote $\xi_t = \left(\xi_{1t},\xi_{2t}^T\right)^T$, where
\begin{flalign*}
\xi_{1t} =\frac{1}{2\sigma_{11}f(0)} \frac{y_{t-1}\mathrm{sign}(\eta_t) }{\sqrt{\omega_0 + \alpha_0 y_{t-1}^2}}, \quad
\xi_{2t} =\Sigma^{-1}\Big( \frac{2y_{t-1}^2}{\omega_0+\alpha_0 y_{t-1}^2} , 2 \Big)^T\big(|\eta_t|-1\big),
\end{flalign*}
and  $\nu = (\nu_1,\nu_2)^T$, where
\begin{flalign*}
\nu_1 =  E\frac{ \omega_0}{2\sqrt{\alpha_0}\left(\omega_0 + \alpha_0 y_{t-1}^2\right)} ,
\quad
\nu_2 =  - E \frac{\sqrt{\alpha_0}}{2\left(\omega_0 + \alpha_0 y_{t-1}^2\right)}.
\end{flalign*}
In the case of $\gamma_0 > 0$, we denote
\begin{flalign*}
\tilde{\xi}_{1t} = \frac{\sqrt{\alpha_0}}{2f(0)}\mathrm{sign}(\eta_t).
\end{flalign*}
The following theorem gives the asymptotic properties of $\hat{\gamma}_n$.
\begin{theorem}\label{theorem4}
Suppose that Assumptions \ref{asm1}--\ref{asm4} hold. \\
$(\mathrm{i})$. If $\gamma_0 < 0$, then
\begin{flalign*}
\sqrt{n}\left(\hat{\gamma}_n - \gamma_0 \right) &= \frac{1}{\sqrt{n}} \sum_{t=1}^n \big(
\zeta_{t} + \mu_1 \xi_{1t} + \mu_2 \nu^T \xi_{2t}\big) + o_p(1)\\
&\stackrel{\mathcal{L}}{\longrightarrow} N\big( 0,~ E\big(\zeta_{t} + \mu_1 \xi_{1t} + \mu_2 \nu^T \xi_{2t}\big)^2 \big).
\end{flalign*}
$(\mathrm{ii})$. If $\gamma_0>0$, then
\begin{flalign*}
\sqrt{n}\left(\hat{\gamma}_n - \gamma_0 \right)&= \frac{1}{\sqrt{n}} \sum_{t=1}^n \big(
\zeta_{t} + \mu_1 \tilde{\xi}_{1t}\big) + o_p(1)\\
&\stackrel{\mathcal{L}}{\longrightarrow}~N\big( 0, ~ E\big(\zeta_{t} + \mu_1 \tilde{\xi}_{1t}\big)^2 \big).
\end{flalign*}
\end{theorem}

Here we provide a heuristic sketch so as to understand why its associated theory involves $f(\cdot)$ and $f'(\cdot)$ and is totally different from that in \cite{fz12}. In the Appendix, we will give detailed proofs.
For illustration, we consider the following toy problem. Suppose that we observe i.i.d. data $\{(X_i,Z_i)^T\}_{i=1}^n$ sampled from  a bivariate normal variable $(X,Z)^T$ with mean $(0,a)^{T}$ and an identity covariance matrix, where $a$ is an unknown parameter and $a \neq 0$. We aim to estimate the quantity $\gamma_1 = E\log|X+a|$, similar to our $\gamma_0$ in a DAR(1) model. Clearly, a natural estimator is 
$$
\hat{\gamma}_1 = \frac{1}{n}\sum_{i=1}^n\log |X_i + \hat{a}_n|,
$$
where $\hat{a}_n = n^{-1}\sum_{i=1}^n Z_i.$ Decompose $\hat{\gamma}_1 - \gamma_1$ as
\begin{flalign*}
\hat{\gamma}_1 - \gamma_1 = \frac{1}{n}\sum_{i=1}^n (\log |X_i + \hat{a}_n|  - \log |X_i + a|) + \frac{1}{n}\sum_{i=1}^n (\log |X_i + a| - \gamma_1): = I_1 + I_2.
\end{flalign*}
A challenge is from $I_1$. Noting that $\log|x+a|=\tfrac{1}{2}\log(x+a)^2$, by Taylor's expansion, the term $I_1$ becomes
\begin{flalign*}
I_1 = \frac{1}{n}\sum_{i=1}^n \frac{1}{X_i+a^*}\, (\hat{a}_n - a).
\end{flalign*}
where $a^*$ is some value between $\hat{a}_n$ and $a$. However, for any constant $c>0$, $n^{-1}\sum_{i=1}^n(X_i+c)^{-1}$ does not converge as $n\rightarrow\infty$. In other words, the Taylor's expansion approach fails.

To solve the mathematical problem for $I_1$, a possible remedy is to apply modern empirical process theory. Denote empirical and true distributions of $X$ by $P_n(x) = n^{-1}\sum_{i=1}^n I(X_i \le x)$
and $P(x) = P(X\le x)$, where $I(A)$ is an indicator function of a set $A$.
Then $I_1 = \int_{x\in R} \log |x|d(P_n(x-\hat{a}_n) - P_n(x-a))$.
Empirical process approximations enable us to obtain that
$$
\sup_{x \in R,\, |a_1-a|\le\delta_n}|P_n(x-a_1)-P(x-a_1)-P_n(x-a)+P(x-a)| = o_p(n^{-1/2-\epsilon}),
$$
provided that $\delta_n \downarrow 0$ at some rate with some $\epsilon>0$, and hence, intuitively,
\begin{eqnarray}\label{approx}
I_1 = \int_{x\in R}\log |x|d(P(x-\hat{a}_n)-P(x-a))+o_p(n^{-1/2}).
\end{eqnarray}
 (Of course, this assertion is not easy to prove and needs to be analyzed carefully since $\log|x|$ is unbounded on $R$. For instance, we added a truncation into our estimator. One reason is to prove (\ref{approx}) in a much easier way.) Through this approximation, $I_1$ reduces to
\begin{flalign*}
I_1 =-\int \log |x| f'(x-a)dx (\hat{a}_n - a) +o_p(n^{-1/2}),
\end{flalign*}
where $f(\cdot)$ is the density of $X$. Define $c = -\int \log |x+a| f'(x)dx$. Then,
\begin{flalign}\label{gamma1}
\sqrt{n} (\hat{\gamma}_1 - \gamma_1) = \frac{1}{\sqrt{n}}\sum_{i=1}^n \{ c (Z_i - a) + \log |X_i + a| - \gamma_1\} + o_p(1).
\end{flalign}
This explains  why  asymptotic variance of $\hat{\gamma}_1$ involves the derivative function $f'(\cdot)$.

 We now discuss how to estimate the asymptotic variance of $ \sqrt{n}(\hat{\gamma}_{n} - \gamma_0)$.  Since it involves the density function $f(\cdot)$ and its derivative $f'(\cdot)$,  the asymptotic expression would not be applicable directly. Similar to Section \ref{twosect}, we propose a random weighting approach for variance estimation. In particular, let $\varpi_1,\cdots,\varpi_n$ be a sequence of i.i.d. nonnegative random variables, with both mean and variance
equal to one.  The definition of $L^*_n(\theta)$ and $\hat{\theta}^*_n$ can be found in (\ref{ladd}). For $t =1,\cdots,n$ and  $k=1,2$, denote
\begin{eqnarray}\label{gammaboots}
\begin{split}
\hat{\eta}^*_t &= \eta_t(\hat{\theta}_n^*),\quad
\mathcal{A}_k^* = \{t: \hat{\eta}_t^* \sqrt{\hat{\alpha}^*_n} - (-1)^k \hat{\phi}_n^*   \in \mathcal{I}_n\}, \\
\hat{\gamma}^*_{nk} &= \frac{1}{\sum_{t \in \mathcal{A}_{k}^*} \varpi_t}
 \sum_{t \in \mathcal{A}_k^*} \varpi_t \log | \hat{\eta}_t^* \sqrt{\hat{\alpha}^*_n} - (-1)^k \hat{\phi}_n^*|.
\end{split}
\end{eqnarray}
The final resampling estimator of $\hat{\gamma}_n$ is defined as
$$
\hat{\gamma}_n^* =(\hat{\gamma}_{n1}^* + \hat{\gamma}_{n2}^*)/2.
$$
To construct the resampling estimator, we insert the random weights into two locations in (\ref{gammaboots}). They correspond to two different parts in (\ref{gamma1}). The following theorem shows that the distribution of $\sqrt{n}(\hat{\gamma}_n - \gamma_0)$ can be approximated by the resampling distribution of $\sqrt{n}(\hat{\gamma}_n^* - \hat{\gamma}_n)$.
\begin{theorem}\label{theorem5}
Suppose that Assumptions \ref{asm1}--\ref{asm4} hold and $\gamma_0 \neq 0$.
Then, conditionally on the data $\{y_0,y_1,\cdots,y_n\}$, in probability,
$\sqrt{n}(\hat{\gamma}_n^* - \hat{\gamma}_n)$
is asymptotically normal as $n \to \infty$,
and its asymptotic variance is the same as that of
$\sqrt{n}(\hat{\gamma}_n - \gamma_0)$.
\end{theorem}

In applications, the asymptotic variance of $\sqrt{n}(\hat{\gamma}_n - \gamma_0)$ can be estimated by using random weighting as follows. First, i.i.d. $\{\varpi_1,\cdots,\varpi_n\}$ are
generated $B$ times from standard exponential distribution, where $B$ is a large number. Each time,  minimize $L_n^*(\theta)$ to obtain $\hat{\theta}^*_n$ and then $\hat{\gamma}_n^*$ is accordingly computed. Denote them as $\hat{\gamma}^{*1}_n,\cdots,\hat{\gamma}^{*B}_n$.
The asymptotic variance $\sigma_{\gamma}^2$ of $\sqrt{n}(\hat{\gamma}_n - \gamma_0)$ is approximated by
\begin{eqnarray*}
\hat{\sigma}_{\gamma}^2 = \frac{1}{B-1}\sum_{b = 1}^B \big(\hat{\gamma}^{*b}_n - \bar{\gamma}^*_n \big)^2 \ \ \mbox{with}
\ \ \bar{\gamma}^*_n = \frac{1}{B}\sum_{b=1}^B \hat{\gamma}^{*b}_n.
\end{eqnarray*}
Once variance estimation is given, confidence intervals can be constructed, as stated in the following corollary.
\begin{corollary}\label{corollary2}
Under the conditions in Theorem \ref{theorem4} and $B = O(n)$, $\hat{\sigma}_{\gamma}^2 \to \sigma_{\gamma}^2$ in probability as $n \to \infty$. Therefore, at the significance level $\underline{\alpha}\in(0, 1)$,
a confidence interval for $\gamma_0$ is
\begin{flalign*}
\left[\hat{\gamma}_n-\frac{\hat{\sigma}_\gamma}{\sqrt{n}}\Phi^{-1}\left(1-\frac{\underline{\alpha}}{2}\right),
~\hat{\gamma}_n+\frac{\hat{\sigma}_\gamma}{\sqrt{n}}\Phi^{-1}\left(1- \frac{\underline{\alpha}}{2}\right )\right],
\end{flalign*}
where  $\Phi(\cdot)$ is  the standard normal distribution.
\end{corollary}

 \subsection{Strict Stationarity Testing}\label{subsecttwo}
Consider the strict stationarity testing
\begin{flalign}\label{testtwo}
H_0:~\gamma_0<0\quad\mbox{v.s.}\quad H_1:~\gamma_0  \ge 0,
\end{flalign}
and
\begin{flalign}\label{testthree}
H_0:~\gamma_0>0\quad\mbox{v.s.}\quad H_1:~\gamma_0 \le 0.
\end{flalign}
With the aid of estimation of $\gamma_0$ and its asymptotics, the above tests are tractable. We give asymptotic critical regions for both testing problems as follows.
\begin{theorem}\label{theorem6}
Suppose that Assumptions \ref{asm1}--\ref{asm4} hold. Let
\begin{eqnarray}
\label{teststat}
T_n=\sqrt{n}\,\frac{\hat{\gamma}_n}{\hat{\sigma}_{\gamma}}
\end{eqnarray}
be the test statistics for (\ref{testtwo}) and (\ref{testthree}), where $\hat{\sigma}_{\gamma}$ is estimated with the resampling size $B = O(n)$. \par
$\mathrm{(i)}$. For the test (\ref{testtwo}), the test defined  by the stationary (ST)
critical region
\begin{flalign*}
\mathrm{C}^{\mathrm{ST}}=\left\{T_n >\Phi^{-1}(1-\underline{\alpha})\right\}
\end{flalign*}
has its asymptotic significance level $\underline{\alpha}$ and is consistent for all $\gamma_0>0$.

$\mathrm{(ii)}$. For the test (\ref{testthree}), the test defined by the nonstationary (NS) critical region
\begin{flalign*}
\mathrm{C}^{\mathrm{NS}}=\{T_n<\Phi^{-1}(\underline{\alpha})\}
\end{flalign*}
has its asymptotic significance level $\underline{\alpha}$ and is consistent for all $\gamma_0<0$.
\end{theorem}

Here we use the normal approximation to construct critical regions for (\ref{testtwo}) and (\ref{testthree}). An alternative approach is to use the resampling distribution of $\hat{\gamma}_n$. If such an approach is applied, the resampling size would need to be larger.

\section{Numerical Studies}\label{foursect}

\subsection{Simulation Studies}\label{4.1sim}
In this subsection, we  conduct numerical studies to assess the finite-sample performance of our GLADE and strict stationarity tests.
In particular, we compare our proposed estimation method with the QMLE studied
in \cite{ling04} and \cite{cll}.

First, we are interested in the performance of the estimators $(\hat{\phi}_n,\hat{\alpha}_n)$ and $\hat{\gamma}_n$ in finite samples. We generate the innovations $\{\eta_t\}$ in the following three scenarios: (a) {\it
the normal distribution $N(0, \pi/2)$}; (b) {\it the Laplace distribution  with density
$f(x)=0.5\exp(-|x|)$};  and (c) {\it the standardized Student's $t_3$ $(\mathrm{st}_3)$ distribution with density $f(x)=4\pi^2/(\pi^2+4x^2)^2$}.
The true parameters are set to be $(\phi_0, \alpha_0, \omega_0)=(0.7, 0.4, 0.5)$ and $(1.0, 3.0, 0.5)$, corresponding to the stationary and explosive cases, respectively. The sample size $n$ is 200 and 400, and 1000 replications are used for each configuration. We take the resampling size $B = 500$ for variance estimation for each configuration. For all cases, the true top Lyapunov exponents have been calculated in Table \ref{table1}.
\begin{table}[!htbp]\scriptsize
\caption{\label{table1} Summary statistics for our proposed procedure under various scenarios based on 1000 replications}
\begin{center}
\begin{tabular}{c|l|cccc|cccc}
   \hline \hline
             &                    & \multicolumn{4}{c|}{$n = 200$}  & \multicolumn{4}{c}{$n=400$}   \\
 \hline
  $\eta$  & Parameters      &  Bias & SE & SEE & CP         &  Bias & SE & SEE & CP \\
 \hline \hline
             & $\phi_0 = 0.700$    & -0.003 & 0.102 & 0.106 & 0.944 & -0.004 & 0.072 & 0.074 & 0.948 \\
Normal       & $\alpha_0 = 0.400$  & -0.014 & 0.085 & 0.082 & 0.925 & -0.007 & 0.056 & 0.058 & 0.951 \\
             & $\omega_0 = 0.500$  & 0.010 & 0.103 & 0.097 & 0.927 & 0.006 & 0.068 & 0.067 & 0.948 \\
             & $\gamma_0 = -0.523$ & -0.005 & 0.109 & 0.113 & 0.958 & -0.005 & 0.074 & 0.076 & 0.950 \\
\hline
             & $\phi_0 = 0.700$     & -0.002 & 0.068 & 0.073 & 0.963 & -0.001 & 0.048 & 0.049 & 0.954 \\
Laplace      & $\alpha_0 = 0.400$   & -0.009 & 0.107 & 0.102 & 0.915 & -0.005 & 0.072 & 0.072 & 0.941 \\
             & $\omega_0 = 0.500$   & 0.019 & 0.135 & 0.129 & 0.930 & 0.013 & 0.092 & 0.090 & 0.947 \\ 
             & $\gamma_0 = -0.440$ & -0.005 & 0.086 & 0.093 & 0.969 & -0.002 & 0.060 & 0.062 & 0.960 \\
\hline
             & $\phi_0 = 0.700$    & -0.002 & 0.082 & 0.085 & 0.943 & 0.001 & 0.057 & 0.058 & 0.940 \\
$\st_3$       & $\alpha_0 = 0.400$  & -0.008 & 0.140 & 0.124 & 0.895 & -0.009 & 0.087 & 0.084 & 0.903 \\
             & $\omega_0 = 0.500$   & 0.012 & 0.147 & 0.139 & 0.906 & 0.016 & 0.111 & 0.103 & 0.919 \\
             & $\gamma_0 = -0.473$ & -0.003 & 0.106 & 0.110 & 0.952 & -0.001 & 0.071 & 0.075 & 0.954 \\
\hline \hline
             & $\phi_0 = 1.000$    & 0.011 & 0.193 & 0.200 & 0.939 & 0.003 & 0.139 & 0.139 & 0.942 \\
Normal       & $\alpha_0 = 3.000$  & -0.009 & 0.338 & 0.326 & 0.932 & -0.006 & 0.223 & 0.228 & 0.952 \\
             & $\omega_0 = 0.500$  & 0.586 & 4.352 & 8.903 & 0.895 & 0.537 & 1.962 & 4.050 & 0.885 \\ 
             & $\gamma_0 = 0.242$  & 0.004 & 0.068 & 0.070 & 0.951 & 0.003 & 0.045 & 0.048 & 0.953 \\
\hline
             & $\phi_0 = 1.000$    & 0.002 & 0.128 & 0.137 & 0.953 & -0.001 & 0.090 & 0.094 & 0.959 \\
Laplace      & $\alpha_0 = 3.000$  & -0.001 & 0.438 & 0.430 & 0.923 & -0.002 & 0.292 & 0.302 & 0.946 \\
             & $\omega_0 = 0.500$  & 1.401 & 7.593 & 80.96 & 0.863 & 2.207 & 25.28 & 28.69 & 0.865 \\  
             & $\gamma_0 = 0.227$  & 0.001 & 0.073 & 0.075 & 0.948 & -0.001 & 0.049 & 0.052 & 0.958 \\
\hline
             & $\phi_0 = 1.000$    & -0.006 & 0.158 & 0.159 & 0.942 & 0.001 & 0.109 & 0.110 & 0.949 \\
$\st_3$       & $\alpha_0 = 3.000$  & 0.022 & 0.536 & 0.521 & 0.921 & 0.001 & 0.375 & 0.356 & 0.926 \\
              & $\omega_0 = 0.500$ & 2.004 & 15.60 & 252.8 & 0.859 & 8.122 & 180.7 & 68.95 & 0.845 \\   
             & $\gamma_0 = 0.183$  & 0.005 & 0.075 & 0.079 & 0.954 & 0.001 & 0.054 & 0.055 & 0.947 \\
\hline \hline
\end{tabular}
\end{center}
Note: Bias and SE are the finite sample bias and standard error of the parameter estimator, SEE means the mean of standard error estimator, and CP is the empirical coverage probabilities of the 95\% confidence intervals.
\end{table}

In Table \ref{table1}, we report the finite-sample biases (Bias), the standard errors (SE), the sample mean of the standard error estimators (SEE) and the empirical coverage probabilities (CP) of the 95\% confidence intervals via normal approximation. Except for $\hat{\omega}_n$ in nonstationary cases, we observe that biases of all other estimators are very small, the variance estimators accurately reflect the true variations and the coverage probabilities agree with the nominal level $95\%$ for almost all cases. These findings confirm that the stationarity assumption is not necessary for estimation of these parameters. They also illustrate that the resampling approach works well for variance estimation. We would like to point out that in some cases, the coverage rates of $\alpha_0$ is slightly lower than the nominal level $95\%$. One possible reason is that $\sqrt{n}(\hat{\alpha}_n - \alpha_0)$ is skewed relative to the normal distribution, especially for small positive $\alpha_0$, since all $\hat{\alpha}_n$ are positive. The basic idea for improvement is to incorporate the skewness information into confidence interval construction. This may be done, for example, by better confidence intervals proposed by \cite{efron}. From  Table \ref{table1}, we can also see that $\hat{\omega}_n$ is not estimable in nonstationary cases.

Next, we evaluate the finite-sample performance of our estimators $(\hat{\phi}_n,\hat{\alpha}_n)$ by comparing with the QMLE. We consider $n=200$ and generate 1000 replications for each configuration. Since both estimation procedures require different conditions on $\eta$ for identification in model (\ref{dar}), it is not a good idea to compare them directly. To make meaningful comparison, we use the average absolute errors (\textsc{AAE}), defined as
\begin{flalign*}
\textsc{AAE}=\frac{1}{2}\Big(\Big|\hat{\phi}-\phi_0\Big|+\Big|\frac{\hat{\alpha}}{\alpha_0}-1\Big|\Big),
\end{flalign*}
where $(\hat{\phi}, \hat{\alpha})$ are the GLADE or QMLE of $(\phi_0, \alpha_0)$, respectively.
Unlike \cite{zhu}, the advantage of this \textsc{AAE} is that it is not necessary to rescale the estimators for both procedures.

\begin{figure}[!hbpt]
\begin{center}
\includegraphics[height= 65mm,width=65mm]{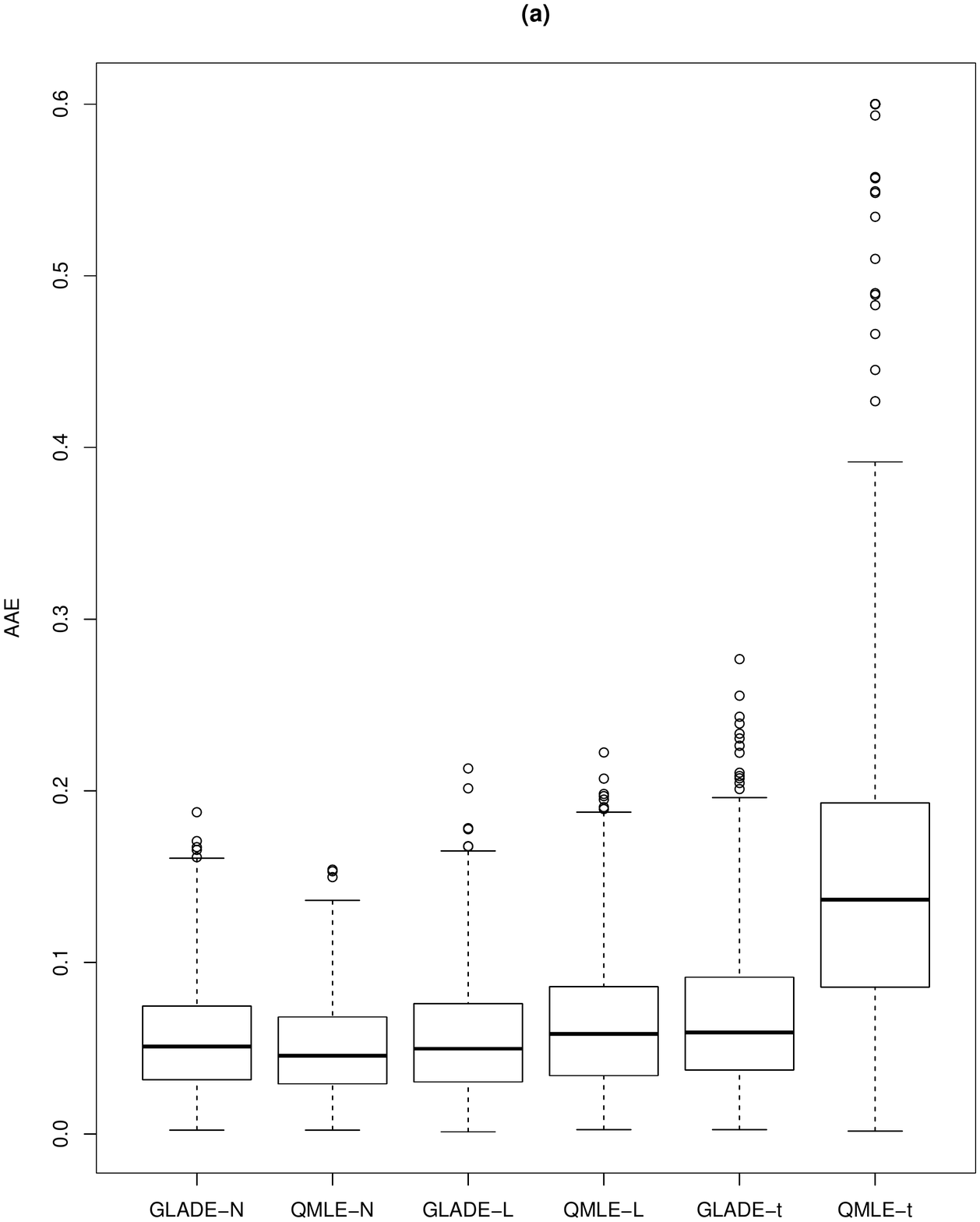}
\includegraphics[height= 65mm,width=65mm]{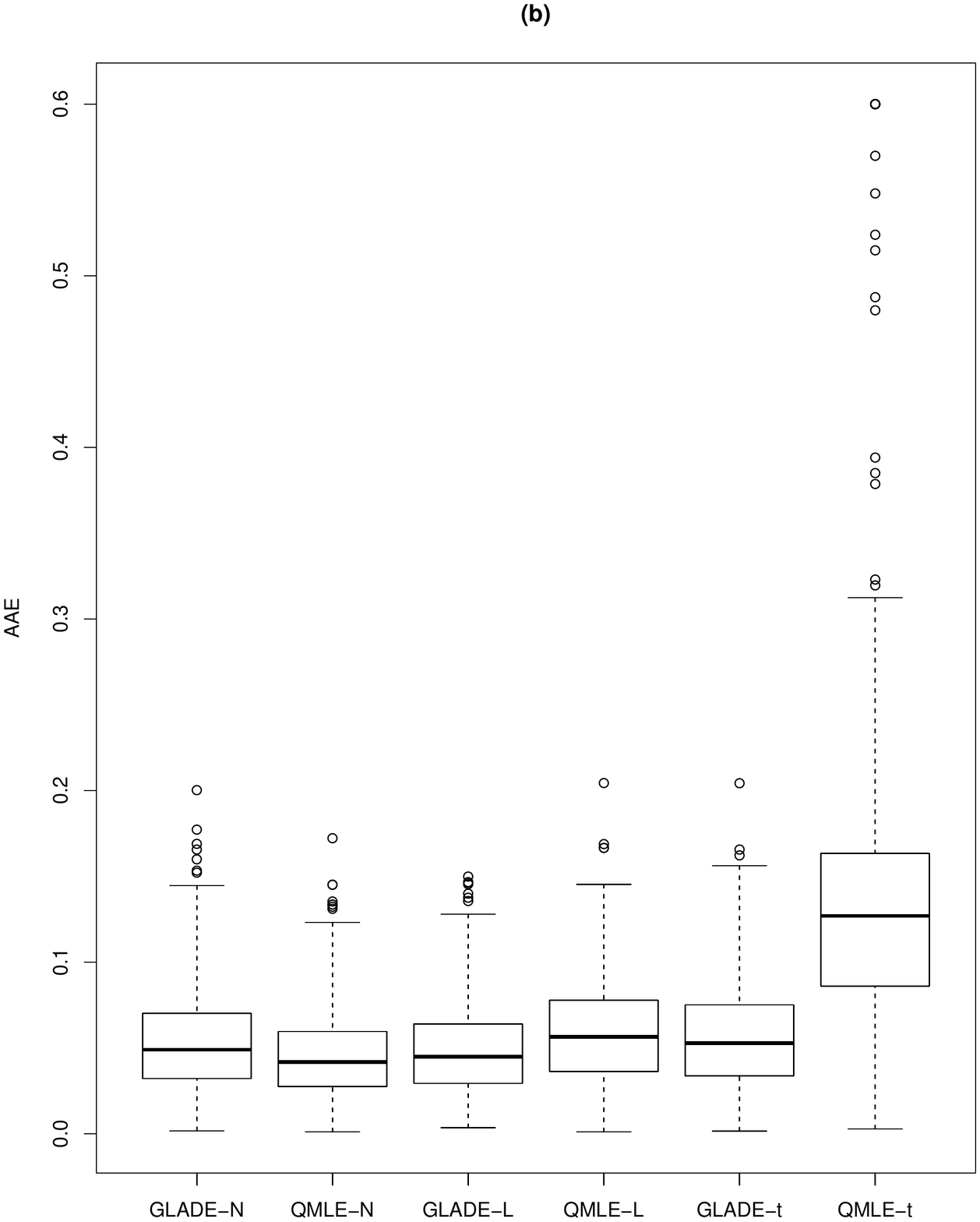}
\caption{Boxplots of the AAE for the GLADE and the QMLE based on 1000 replications. (a) $(\phi_0, \alpha_0, \omega_0)=
(0.7, 0.4, 0.5)$ and (b) $(\phi_0, \alpha_0, \omega_0)=(1.0, 3.0, 0.5)$. `GLADE-N',`GLADE-L' and `GLADE-t' mean `GLADE's when errors are normal, Laplace and $\mathrm{st}_3$, respectively. `QMLE-N',QMLE-L' and `QMLE-t' are defined similarly.}\label{fig2}
\end{center}
\end{figure}

Fig.\ref{fig2} shows the box-plots of the AAE for the GLADE and QMLE.
It is observed from Fig. \ref{fig2} that the GLADE outperforms the QMLE when the innovation $\eta$ follows  $\mathrm{st}_3$ and Laplace distribution, thereby supporting the robust properties of our GLADE, either stationary or explosive. In particular, the QMLE performs much worse when errors have infinite fourth moments. We also observe that the QMLE with normal errors behaves better. It is not surprising since in this case, the QMLE is efficient, while the GLADE is not. 

Finally, we illustrate the performance of strict stationarity tests developed in Section \ref{threesect}. We keep the standardized Student's $t_3$ distribution for $\eta_t$, $\omega_0 = 0.5$ but $\phi_0$ varies from 0.6 to 1.3 and $\alpha_0 = 2\phi_0$. In this scenario, we have $\gamma_0 =0$ for $(\phi_0,\alpha_0)= (0.922,1.844)$ and hence, $\gamma_0 >0$ or $<0 $ if $\phi_0 > 0.922$ or $< 0.922$, respectively. We take the sample size $n=200, 400, 800$, and 1000 replications for each configuration.

Tables \ref{table2} and \ref{table3} summarize empirical frequencies of rejection for (\ref{testtwo}) and (\ref{testthree}) for various values of $\phi_0$, respectively. 
It is observed that when $\gamma_0 =0$, the rejection frequencies of the two tests agree with the nominal level $5\%$ as $n$ increases. As expected, the frequency of rejection of the $\mathrm{C^{ST}}$ test increases with $\gamma_0$, while that of the $\mathrm{C^{NS}}$ test decreases. Overall, the power of the two tests is significant, as shown in theory.
\begin{table}[!htbp]
\begin{center}
\caption{\label{table2} Relative frequency of rejection of the test (\ref{testtwo}): $H_0:  \gamma_0 < 0$ with $(\omega_0,\alpha_0)= (0.5, 2\phi_0)$ based on 1000 replications.}
\begin{tabular}{ccccccccc}
\hline \hline
 &  \multicolumn{7}{c}{$\phi_0$}  \\
 \hline
 &  0.600 &0.700&0.800&0.922&1.000&1.100&1.300 \\
\hline
$n=200$  & 0.000 & 0.000 & 0.001 & 0.053 & 0.173 & 0.528 & 0.943 \\
$n=400$  & 0.000 & 0.000 & 0.000 & 0.062 & 0.252 & 0.754 & 1.000 \\
$n=800$  & 0.000 & 0.000 & 0.000 & 0.057 & 0.452 & 0.957 & 1.000 \\
\hline \hline
\end{tabular}
\end{center}
\end{table}
\begin{table}[!htbp]
\centering
\caption{ \label{table3} Relative frequency of rejection of the test (\ref{testthree}): $H_0:\gamma_0 > 0$  with $(\omega_0,\alpha_0)=(0.5, 2\phi_0)$ based on 1000 replications.}
\begin{tabular}{ccccccccc}
\hline \hline
 &  \multicolumn{7}{c}{$\phi_0$}  \\
 \hline
 & 0.600 &0.700&0.800&0.922&1.000&1.100&1.300 \\
\hline
$n=200$  & 0.954 & 0.713 & 0.313 & 0.038 & 0.009 & 0.000 & 0.000 \\
$n=400$  & 0.996 & 0.933 & 0.501 & 0.048 & 0.003 & 0.000 & 0.000 \\
$n=800$  & 1.000 & 1.000 & 0.796 & 0.054 & 0.002 & 0.000 & 0.000 \\
\hline \hline
\end{tabular}
\end{table}

\subsection{An Empirical Study}
\begin{figure}[!htbp]
 \begin{center}
\includegraphics[height=8cm,width=13cm]{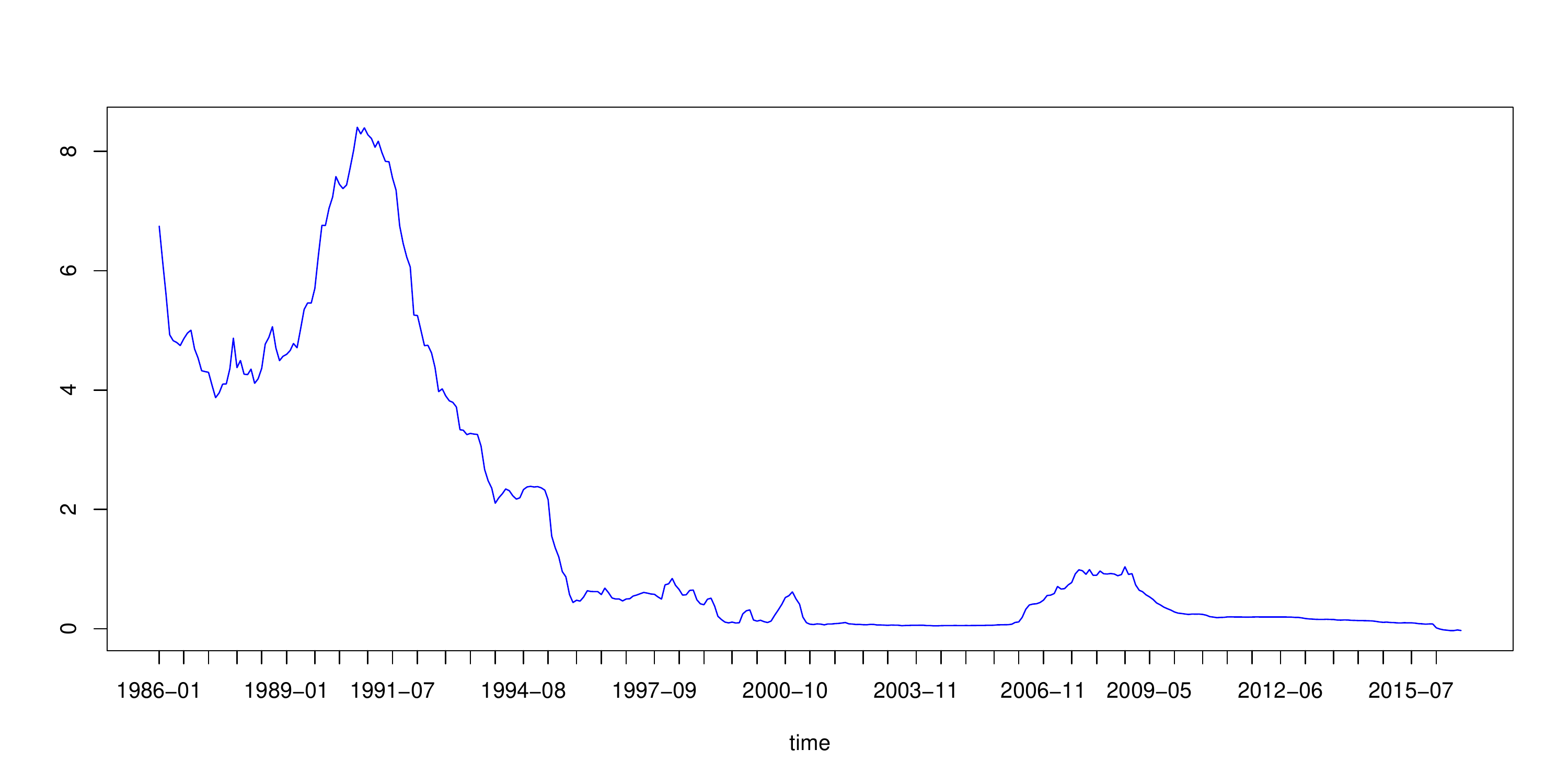}
\end{center}
\caption{The monthly 3-Month LIBOR (in percent), based on Japanese Yen, from Jan 1986 to Sep 2016.}\label{fig3}
\end{figure}

To illustrate the use of our proposed inference procedure, we consider a data set of the monthly 3-Month London Interbank Offered Rate (LIBOR) $\{y_t\}$ (based on Japanese Yen) over the period from January 1986 to September 2016. The observed level series of LIBOR is plotted in Fig. \ref{fig3}. 
Here we aim to check whether this level series over the given period is
strictly stationary by fitting a DAR(1) model. If the significance level is set to be 5\%,
then the stationarity hypothesis of test (\ref{testtwo}) can be tested by comparing
the $t$-type statistic $T_n$ with 1.64 and rejecting the null hypothesis if $T_n> 1.64$.
To calculate $T_n$, we set the resampling size to be 1000 for variance estimation of $\hat{\gamma}_n$.

We fit the level series directly instead of its log-return series since this data set contains negative observations around year 2016. Table \ref{table4} reports estimates of $\theta$, the top Lyaponov exponent $\gamma$ and the value of test statistics $T_n$ for (\ref{testtwo}).
\begin{table}[!htbp]
\centering
\caption{\label{table4} Parameter estimates and test statistic $T_n$ of the stationarity test in (\ref{testtwo}).}
\begin{tabular}{c|ccccc}
\hline \hline
Index & $\hat{\phi}$ &$\hat{\alpha} $&$\hat{\omega}$ & $\hat{\gamma}$ & $T_n$ \\
LIBOR & 0.9940&0.0086&$1.5454\times10^{-5}$ &-0.0237 &-4.7760\\
&(0.0036)&0.0026&($0.7976\times 10^{-5}$)&(0.0050)&(1.0000)\\
\hline \hline
\end{tabular}\\
{Note: the values in the first four parenthesises are the standard deviations of parameter estimates and the value in the last parenthesis is the $p$-value of $T_n$.}
\end{table}

From Fig.\,\ref{fig3}, the level LIBOR series seems to be nonstationary intuitively.
Meanwhile, the value of $\hat{\phi}$ in Table \ref{table4} is insignificant from 1. However, we should do the judgement with more caution in the DAR model since a unit root in the conditional mean does not necessarily imply nonstationarity as lagged values of the process enter into the conditional variance. As a matter of fact, the test statistic $T_n$ is $-3.4209$, which is less than 1.64, meaning that we do not have significant evidence to reject the null hypothesis of strict stationarity ($p$-value=0.999).
This example demonstrates that the series constructed from the DAR(1) model is allowed to be persistent but remains strictly stationary, see \cite{NR2014}. It further illustrates that it may not be reliable to judge the stationarity of nonlinear time series processes using the traditional unit-root tests applied to linear AR  models, or by visually inspecting the data.

\begin{figure}[!htbp]
\begin{center}
\includegraphics[height=10cm,width=15cm]{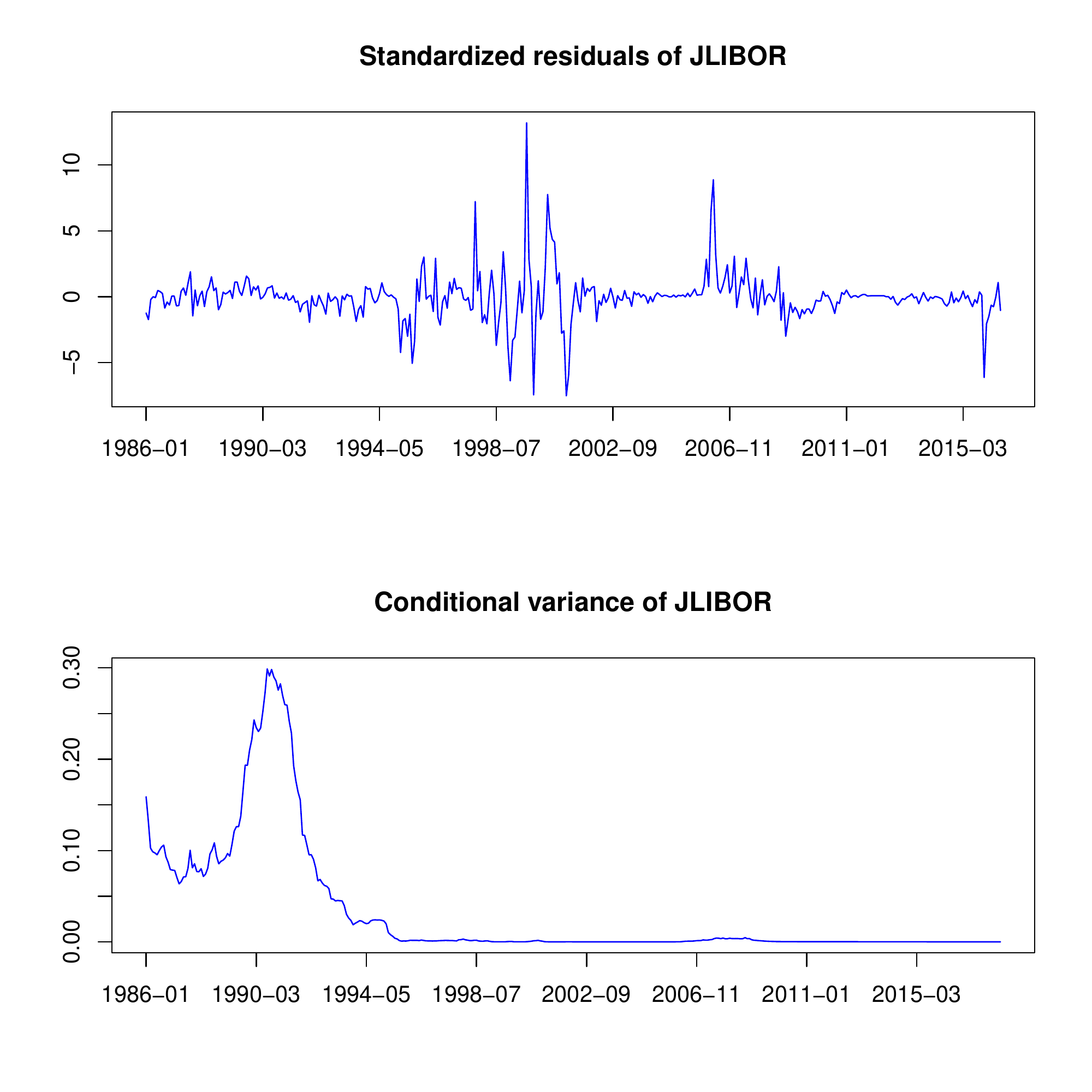}
\end{center}
\caption{The standardized residuals and conditional variance from model (\ref{simulatedmodel}).}
\label{fig4}
\end{figure}

The final estimated model is presented as following:
\begin{flalign}\label{simulatedmodel}
y_t=0.994 y_{t-1}+\eta_t\sqrt{1.5454\times10^{-5}+0.0086y_{t-1}^{2}},\quad t=2,...,369.
\end{flalign}
The standardized residuals and conditional variance from model (\ref{simulatedmodel}) are plotted in Fig. \ref{fig4}.

\begin{figure}[!htbp]
\begin{center}
\includegraphics[height=8cm,width=13cm]{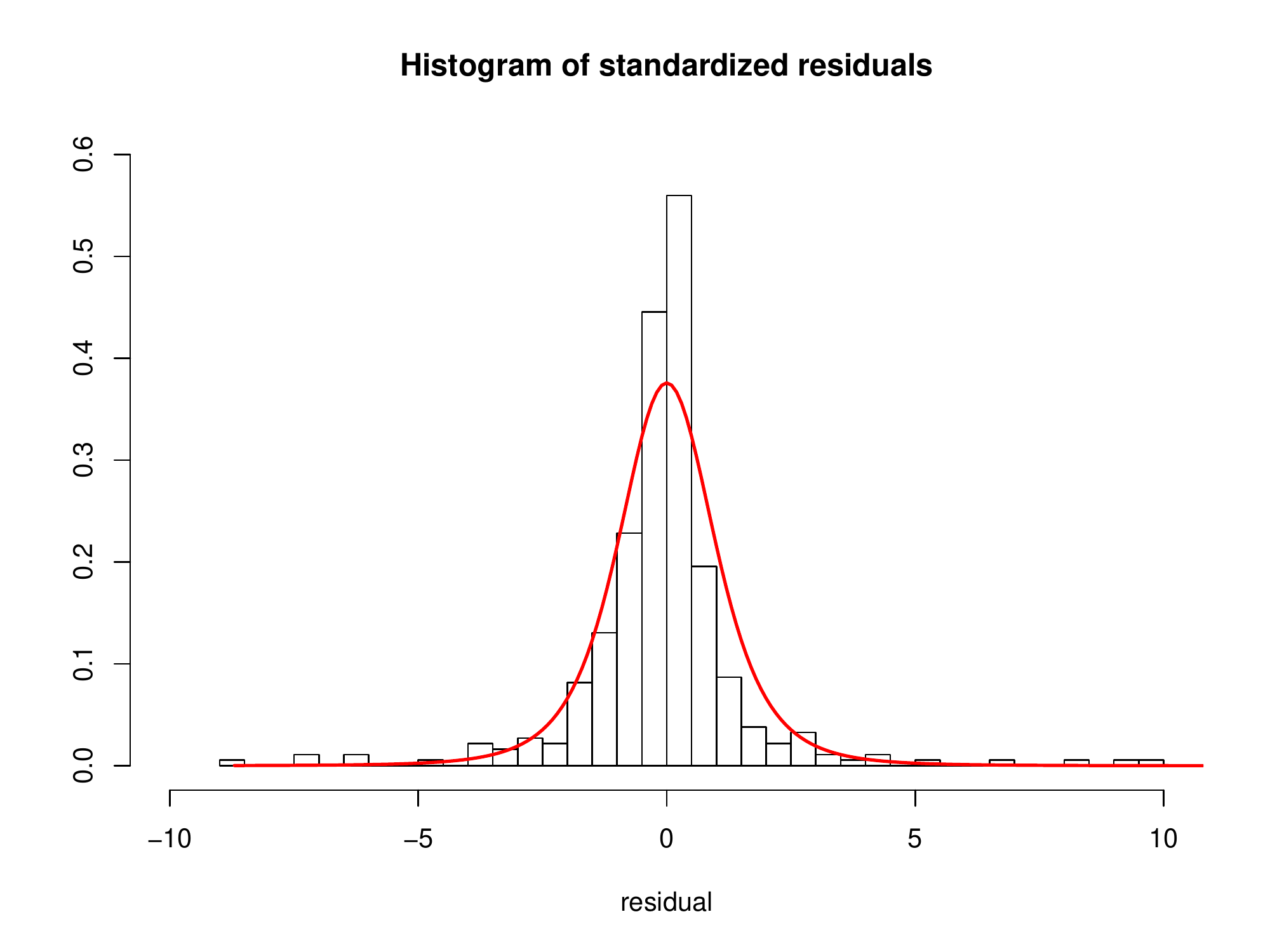}
\end{center}
\caption{The histogram of the residuals $\{\hat{\eta}_t\}$ from model (\ref{simulatedmodel}).} \label{fig5}
\end{figure}

\begin{figure}[!htbp]
 \begin{center}
\includegraphics[height=8cm,width=13cm]{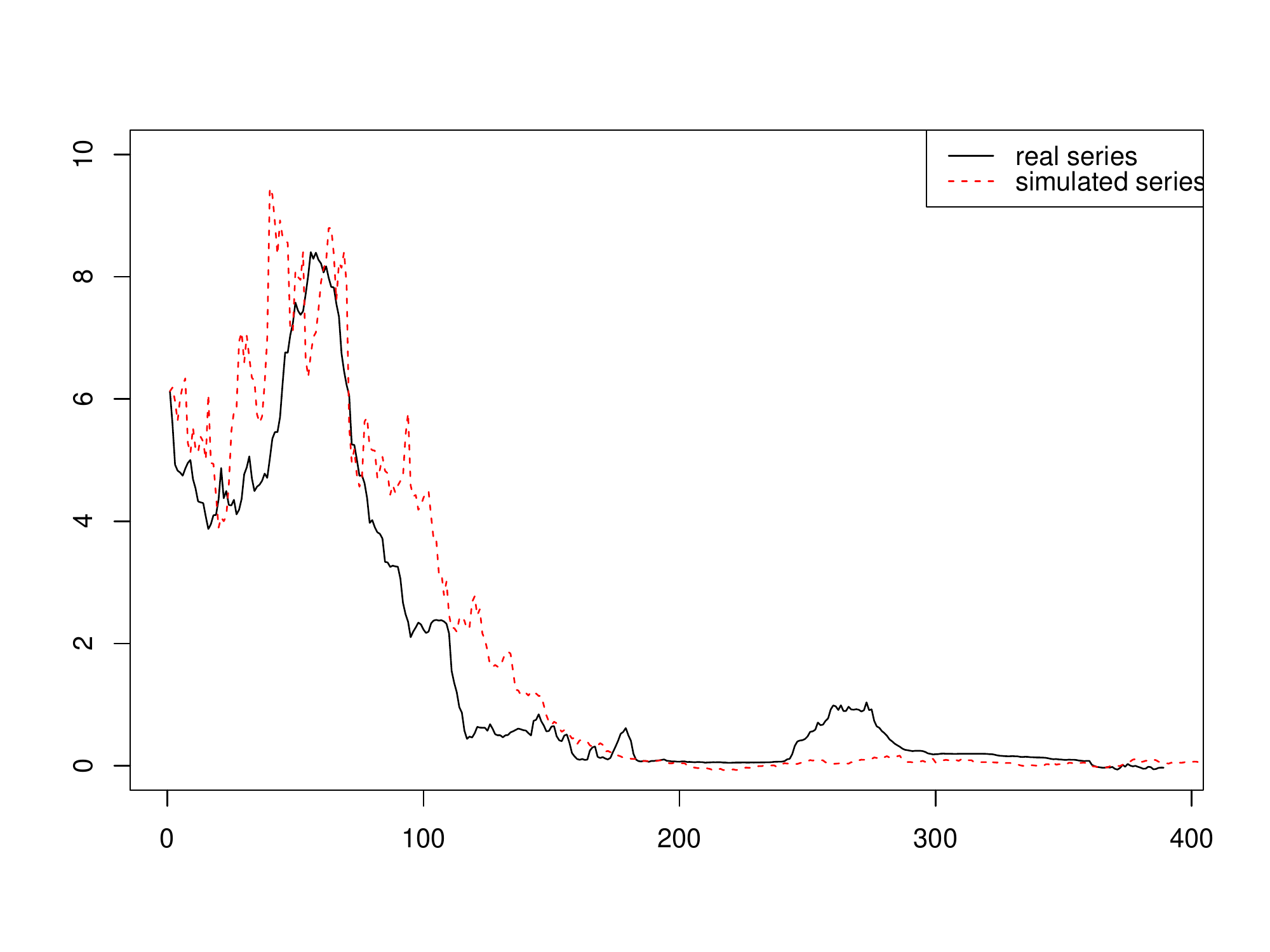}
\end{center}
\caption{The simulated series with $\hat{\phi}$, $\hat{\alpha}$ and $\hat{\omega}$ and
$t_3$ distributed innovations.} \label{fig6}
\end{figure}

To assess the performance of the fitting DAR(1) model, we simulate a sample path from the estimated DAR(1) model (\ref{simulatedmodel}) with an initial value $y_1$, where $\{\eta_t\}$ is i.i.d. Student's $t_3$.
Here, the Student's $t$ distribution is chosen based on the histogram of the residuals in Fig. \ref{fig5} and the degree of freedom 3 is approximately obtained from the standard deviation (i.e., 1.932) of the residuals.
Clearly, the model (\ref{simulatedmodel}) is strictly stationary according to the stationarity region,
see, e.g., \cite{bk} and \cite{cll}. A simulated path is plotted in Fig. \ref{fig6},
from which we can see that the simulated series behaves like that of the true LIBOR level series.
This further provides an evidence for our findings.


\section{Conclusion}\label{fivesect}
Testing strict stationarity is important in the context of nonlinear time series analysis.
This paper addresses the strict stationarity testing  on DAR models and develops a unified framework for the inference of both stationary and explosive DAR(1) processes.
We also propose a global LADE of $(\phi_0, \alpha_0)$ in the DAR(1) model and establish its asymptotic theory without strict stationarity assumption.
If one is interested in the inference on $\phi_0$ or $\alpha_0$, strict stationarity testing is unnecessary.

It is worth noting that one of  assumptions is $\eta$ is symmetric, which is  met by most commonly used errors. 
When this assumption is violated, the statistical inference of DAR models is absent since the dynamic behavior of $\{y_t\}$ is unknown in the nonstationary case.
We leave it to future studies.

\section*{Acknowledgements}
The authors would like to thank the Co-Editor, the Associate Editor and two anonymous reviewers for their helpful comments and valuable suggestions. Shaojun Guo acknowledges partially  the support of  the NSFC (No. 11771447)
and the Fundamental Research Funds for the Central Universities and the Research Funds of Renmin University of China.
Dong Li acknowledges partially the support of the NSFC (No.11571348 and No.11771239).
Muyi Li acknowledges the support of the NSFC (No.71671150), Fujian Key Laboratory of Statistical Science and IT service in both WISE and School of Economics, Xiamen University.

\appendix
\section*{Appendix}
\section{Proofs of Theorems}
Denote by $M$ any positive constant whose
value is unimportant and can be different throughout the proofs.
Recall the fact that
\begin{flalign}\label{fact}
|y_t|/\rho^t\rightarrow\infty\quad \mbox{a.s. as $t\rightarrow\infty$}
\end{flalign}
for any $\rho\in(1, e^{\gamma_0})$  by Theorem 1 in \cite{cll}.
\subsection{Proof of Theorem \ref{consistency}}
The result stated in (i) is standard, see \cite{zhu}. Consider the case (ii). Clearly,
$(\hat{\phi}_n, \hat{\alpha}_n, \hat{\omega}_n)^T=\arg\min_{\theta\in\Theta} Q_n(\theta)$, where
$Q_n(\theta)=\{L_n(\theta)-L_n(\theta_0)\}/n$. We have
\begin{flalign*}
Q_n(\theta)=O_n(\phi,\alpha)+R_{1n}(\theta)+R_{2n}(\theta),
\end{flalign*}
where
\begin{flalign*}
O_n(\phi,\alpha)=\frac{1}{n}\sum_{t=1}^n\left\{\frac{1}{2}\log\frac{\alpha}{\alpha_0}+
\frac{|\eta_t\mathrm{sign}(y_{t-1})\sqrt{\alpha_0}-(\phi-\phi_0)|}{\sqrt{\alpha}}-|\eta_t|\right\}
\end{flalign*}
and
\begin{flalign*}
R_{1n}(\theta)&=\frac{1}{2n}\sum_{t=1}^n\log\frac{\alpha_0(\omega+\alpha y_{t-1}^2)}
{\alpha(\omega_0+\alpha_0 y_{t-1}^2)},\\
R_{2n}(\theta)&=\frac{1}{2n}\sum_{t=1}^n\bigg\{\frac{|\eta_t\sqrt{\omega_0+\alpha_0y^2_{t-1}}
-(\phi-\phi_0)y_{t-1}|}{\sqrt{\omega+\alpha y^2_{t-1}}}
-\frac{|\eta_t\mathrm{sign}(y_{t-1})\sqrt{\alpha_0}-(\phi-\phi_0)|}{\sqrt{\alpha}}\bigg\}.
\end{flalign*}
Note that $(\eta_1\mathrm{sign}(y_{0}),...,\eta_n\mathrm{sign}(y_{n-1}))$ and $(\eta_1,...,\eta_n)$
have the same distribution. Thus,
\begin{flalign*}
\lim_{n\rightarrow\infty}O_n(\phi,\alpha)&=\tfrac{1}{2}\log\frac{\alpha}{\alpha_0}+
\frac{E|\eta_t\sqrt{\alpha_0}-(\phi-\phi_0)|}{\sqrt{\alpha}}-1\quad \quad\mbox{a.s.}\\
&\geq\tfrac{1}{2}\log(\alpha/\alpha_0)+
\sqrt{\alpha_0/\alpha}-1\geq0
\end{flalign*}
because $E|\eta_t-c|\geq E|\eta_t|$ for any $c\in R$ and the inequality $\log x\leq x-1$ for $x>0$.
The equality holds if and only if $\phi=\phi_0$ and $\alpha=\alpha_0$.

For $R_{1n}(\theta)$, by the mean value theorem, the compactness of $\Theta$ and (\ref{fact}), we have
\begin{flalign*}
\lim_{n\rightarrow\infty}\sup_{\theta\in\Theta}|R_{1n}(\theta)|\leq
M\lim_{n\rightarrow\infty}\frac{1}{n}\sum_{t=1}^{n}\frac{1}{1+y^2_{t-1}}\rightarrow0\quad\mbox{a.s.}
\end{flalign*}

For $R_{2n}(\theta)$, by the compactness of $\Theta$ and (\ref{fact}), a simple calculus yields
\begin{flalign*}
\lim_{n\rightarrow\infty}\sup_{\theta\in\Theta}|R_{2n}(\theta)|\leq
M\lim_{n\rightarrow\infty}\frac{1}{n}\sum_{t=1}^{n}\frac{|\eta_t|}{|y_{t-1}|}\rightarrow0\quad\mbox{a.s.}
\end{flalign*}
Thus, the proof can be completed by standard arguments, invoking the compactness of $\Theta$.

\subsection{Proof of Theorem \ref{theorem2}}
(i). When $\gamma<0$, $\{y_t\}$ is strictly stationary and ergodic,
and there exists some $\iota>0$ such that $E|y_t|^\iota<\infty$, see \cite{bk}.
Thus, Assumption 2  in \cite{zhu} is satisfied and in turn the result holds.

(ii). We first reparameterize the objective function as
\begin{flalign}\label{dongobject}
H_n(u)=L_n(\theta_0+u)-L_n(\theta_0),
\end{flalign}
where $u\in\Lambda:=\{u=(u_1, u_2, u_3)^T: u+\theta_0\in \Theta\}$.

Let $\hat{u}_n=\hat{\theta}_n-\theta_0:=(\hat{u}_{1n}, \hat{u}_{2n}, \hat{u}_{3n})^T
=(\hat{v}_n^T, \hat{u}_{3n})^T$.
We can see that $\hat{u}_n$ is the minimizer of $H_n(u)$ in $\Lambda$.
By Theorem \ref{consistency}(ii), $\hat{v}_n=o_p(1)$.
Using the similar arguments as in \cite{zhling,zhu}, by the fact (\ref{fact}), we have
\begin{flalign*}
H_n(\hat{u}_n)=(\sqrt{n}\hat{v}_n)^TT_n+(\sqrt{n}\hat{v}_n)^T\Omega(\sqrt{n}\hat{v}_n)+
o_p(\sqrt{n}\|\hat{v}_n\|+n\|\hat{v}_n\|^2),
\end{flalign*}
where
\begin{flalign*}
T_n=-\frac{1}{\sqrt{n}}\sum_{t=1}^n\Big(\frac{\mathrm{sign}(\eta_t)}{\sqrt{\alpha_0}},
\frac{|\eta_t|-1}{2\alpha_0}\Big)^T\,\,\mbox{and}\,\,
\Omega=\mathrm{diag}(f(0)/\alpha_0,\: 1/(8\alpha_0^2)).
\end{flalign*}
By the central limit theorem, we have $T_n\stackrel{\mathcal{L}}{\longrightarrow} N(0, \mathrm{diag}(1/\alpha_0, \,(\kappa_\eta-1)/(4\alpha_0^2)))$.
Let $\lambda_{\mathrm{min}}=\min\{f(0)/\alpha_0,\: 1/(8\alpha_0^2)\}>0$. Then
\begin{flalign*}
H_n(\hat{u}_n)\geq -\|\sqrt{n}\hat{v}_n\|\{\|T_n\|+o_p(1)\}+\|\sqrt{n}\hat{v}_n\|^2\{\lambda_{\mathrm{min}}+o_p(1)\}.
\end{flalign*}
Note that $H_n(\hat{u}_n)\leq 0$ by the definition of $\hat{\theta}_n$. Thus,
\begin{flalign*}
\|\sqrt{n}\hat{v}_n\|\leq \{\lambda_{\mathrm{min}}+o_p(1)\}^{-1}\{\|T_n\|+o_p(1)\}=O_p(1).
\end{flalign*}
Next, let $v_n^*=-\Omega^{-1}T_n/(2\sqrt{n})$. Then,
\begin{flalign*}
\sqrt{n}v_n^*=-\Omega^{-1}T_n/2\stackrel{\mathcal{L}}{\longrightarrow}N(0, \mathcal{J}_N).
\end{flalign*}
Using the previous facts, a simple calculus gives that
\begin{flalign*}
H_n(\hat{u}_n)-H_n(u_n^*)&=(\sqrt{n}\hat{v}_n-\sqrt{n}v_n^*)^T\Omega(\sqrt{n}\hat{v}_n-\sqrt{n}v_n^*)+o_p(1)\\
&\geq \lambda_{\mathrm{min}}\|\sqrt{n}\hat{v}_n-\sqrt{n}v_n^*\|^2+o_p(1),
\end{flalign*}
where $u_n^*=(v_n^{*T}, 0)^T$. Note that
$H_n(\hat{u}_n)-H_n(u_n^*)=L_n(\theta_0+\hat{u}_n)-L_n(\theta_0+u^*_n)\leq0$ a.s. Thus, we have
$\|\sqrt{n}\hat{v}_n-\sqrt{n}v_n^*\|=o_p(1).$

Finally, we have $\sqrt n(\hat{\phi}_n-\phi_0, ~\hat{\alpha}_n-\alpha_0)^T
=\sqrt{n}\hat{v}_n=\sqrt{n}v_n^*+o_p(1)$
and then (ii) holds.

\subsection{Proof of Theorem \ref{theorem3}}
We first prove (ii).  Let $P^*$  be the joint probability of $(\varpi_1,\cdots, \varpi_n)$ and $(y_1,\cdots,y_n)$. 
Replacing $L_n(\theta)$ in (\ref{dongobject}) by $L_n^*(\theta)$ and repeating the proof of 
Theorem \ref{theorem2}, we have that
\begin{flalign*}
\sqrt{n}( \hat{\phi}_n^* - \phi_0 ) &= \frac{\sqrt{\alpha_0}}{2 f(0)} \frac{1}{\sqrt{n}}
\sum_{t=1}^n \varpi_t\sgn(\eta_t) +  o_{P^*}(1),\\
\sqrt{n}( \hat{\alpha}_n^* - \alpha_0 ) &=
\frac{2 {\alpha_0}}{\sqrt{n}}\sum_{t=1}^n \varpi_t \{|\eta_t|-1\} +  o_{P^*}(1).
\end{flalign*}
For more details, the readers can refer to the recent paper of \cite{Dovonon}.
Note that $( \hat{\phi}_n^* - \hat{\phi}_n) = ( \hat{\phi}_n^* - \phi_0 ) - ( \hat{\phi}_n - \phi_0 )$ and so does for $( \hat{\alpha}_n^* - \hat{\alpha}_n)$. Hence, together with asymptotic representations of $\sqrt{n}( \hat{\phi}_n - \phi_0 )$
 and $ \sqrt{n}( \hat{\alpha}_n - \alpha_0 ) $ in the proof of Theorem \ref{theorem2}, we immediately obtain that
\begin{flalign*}
\sqrt{n}( \hat{\phi}_n^* - \hat{\phi}_n ) &= \frac{\sqrt{\alpha_0}}{2 f(0)} \frac{1}{\sqrt{n}}
\sum_{t=1}^n (\varpi_t - 1)\sgn(\eta_t) +  o_{P^*}(1),\\
\sqrt{n}( \hat{\alpha}_n^* - \hat{\alpha}_n ) &=
\frac{2 {\alpha_0}}{\sqrt{n}}\sum_{t=1}^n (\varpi_t - 1) \{|\eta_t|-1\} +  o_{P^*}(1),
\end{flalign*}
and therefore, conditional on the data $\{y_0,y_1,\cdots,y_n\}$, the part (ii) follows. The proof of (i) is proved analogously.

\subsection{Proof of Theorem \ref{theorem4}}
First, we decompose $\hat{\phi}_n + \hat{\eta}_t \sqrt{\hat{\alpha}_n}$ for each $t =1,\cdots,n$.
In fact, $\hat{\phi}_n + \hat{\eta}_t \sqrt{\hat{\alpha}_n}$ is rewritten as
\begin{flalign*}
\hat{\phi}_n + \hat{\eta}_t \sqrt{\hat{\alpha}_n} &= \hat{\phi}_n + \frac{y_t - \hat{\phi}_n y_{t-1}}{\sqrt {\hat{\omega}_n + \hat{\alpha}_n y_{t-1}^2}}\sqrt{\hat{\alpha}_n} \\
&=\hat{\phi}_n + \frac{(\phi_0 - \hat{\phi}_n) y_{t-1}}{\sqrt {\hat{\omega}_n + \hat{\alpha}_n y_{t-1}^2}}\sqrt{\hat{\alpha}_n} + \sqrt {\frac{{\omega}_0 + {\alpha}_0 y_{t-1}^2}{\hat{\omega}_n + \hat{\alpha}_n y_{t-1}^2}}\eta_t\sqrt{\hat{\alpha}_n}\\
& = \phi_0 + \eta_t \sqrt{\alpha_0} + v_{t}^{(1)}(\hat{\theta}_n) + \eta_t u_{t}(\hat{\theta}_n).
\end{flalign*}
where $\hat{\theta}_n = (\hat{\phi}_n,\hat{\alpha}_n,\hat{\omega}_n)^T$,
\begin{flalign*}
v_{t}^{(1)}(\theta) = (\phi - \phi_0)\Big( 1 - \frac{ y_{t-1}\sqrt{\alpha}}{\sqrt {\omega + \alpha y_{t-1}^2}}\Big),\quad
u_{t}(\theta) =\sqrt {\frac{\omega_0 +\alpha_0 y_{t-1}^2 }{\omega + \alpha y_{t-1}^2}}\sqrt{\alpha} - \sqrt{\alpha_0}.
\end{flalign*}
Similarly,
\begin{eqnarray*}
 \hat{\eta}_t \sqrt{\hat{\alpha}_n} - \hat{\phi}_n
 =  \eta_t \sqrt{\alpha_0}  - \phi_0 + v_{t}^{(2)}(\hat{\theta}_n) + \eta_t u_{t}(\hat{\theta}_n),
\end{eqnarray*}
where $v_{t}^{(2)}(\theta) = -(\phi - \phi_0)\Big( 1 + \frac{ y_{t-1}\sqrt{\alpha}}{\sqrt {\omega + \alpha y_{t-1}^2}}\Big).
$

(I) Consider the stationary case, i.e. $\gamma_0 < 0$.  Let $v_{nt}^{(k)}(\vartheta ) = n^{1/2} v_{t}^{(1)}(\theta_0 + n^{-1/2} \vartheta )$, $k=1,2$, and $u_{nt}(\vartheta ) = n^{1/2} u_{t}(\theta_0 + n^{-1/2} \vartheta )$. For $k =1,2$, define
\begin{eqnarray*}
H_{nk}(x;\vartheta) &=& \frac{1}{n}\sum_{t = 1}^n I\Big( \eta_t  \le \frac{x + (-1)^k \phi_0 - n^{-1/2}v_{nt}^{(k)}(\vartheta)}{\sqrt{\alpha_0} + n^{-1/2}u_{nt}(\vartheta)}\Big), \\
H_k(x;\vartheta) &=& \frac{1}{n}\sum_{t = 1}^n F\Big(\frac{x + (-1)^k\phi_0 - n^{-1/2}v_{nt}^{(k)}(\vartheta)}{\sqrt{\alpha_0} + n^{-1/2} u_{nt}(\vartheta)}\Big).
\end{eqnarray*}
If we denote $\hat{\vartheta}_n = n^{1/2}(\hat{\theta}_n - \theta_0)$ and $\mathcal{I}_n = [-n^2, -n^{-2}]\cup [n^{-2},n^2]$, then
\begin{eqnarray*}
2\hat{\gamma}_n =   \int_{\mathcal{I}_n} \log |x| d H_{n1}(x;\hat{\vartheta}_n) + \int_{\mathcal{I}_n} \log |x| d H_{n2}(x;\hat{\vartheta}_n)
= \hat{\gamma}_{n1} + \hat{\gamma}_{n2}.
\end{eqnarray*}
We decompose $\hat{\gamma}_{nk} - \gamma_0$ into 
\begin{eqnarray*}
\hat{\gamma}_{nk} - \gamma_0 &=&
\int_{\mathcal{I}_n} \log|x| d \big( H_{nk}(x;\hat{\vartheta}_n) - H_{nk}(x;0) - H_k(x;\hat{\vartheta}_n ) + H_k(x;0)\big)\\
&& + \int_{\mathcal{I}_n} \log |x| d(H_{nk}(x;0) - H_k(x;0)) + \int_{R/\mathcal{I}_n} \log |x| d H_k(x;0)  \\
&& + \int_{\mathcal{I}_n} \log |x| d(H_k(x;\hat{\vartheta}_n) - H_k(x;0))\\
&=& I_{n1}^{(k)} + I_{n2}^{(k)} + I_{n3}^{(k)} + I_{n4}^{(k)}.
\end{eqnarray*}
In the following, we will show that
\begin{flalign*}
&(a)\quad \sqrt{n} I_{n1}^{(k)} = o_p(1),\quad k = 1,2;\\
&(b)\quad \hskip 0.1cm
\sqrt{n} \big ( I_{n2}^{(1)} + I_{n2}^{(2)} \big )= \frac{1}{\sqrt{n}} \sum_{t=1}^n
\big (\zeta_{1t} + \zeta_{2t}\big ) + o_p(1); \\
&(c)\quad \sqrt{n} I_{n3}^{(k)} = o_p(1),\quad k = 1,2; \\
&(d)\quad \hskip 0.1cm
\sqrt{n} \big ( I_{n4}^{(1)} + I_{n4}^{(2)} \big ) = \frac{1}{\sqrt{n}} \sum_{t=1}^n
\big ( 2\mu_1 \xi_{1t} + 2 \mu_2 \nu^T \xi_{2t}\big)+ o_p(1). &&
\end{flalign*}

Suppose that (a)-(d) are proved. It immediately follows that
\begin{flalign*}
\sqrt{n}\left(\hat{\gamma}_n - \gamma_0 \right) = \frac{1}{\sqrt{n}} \sum_{t=1}^n \big(
\zeta_{t} + \mu_1 \xi_{1t} + \mu_2 \nu^T \xi_{2t}\big) + o_p(1),
\end{flalign*}
and hence $\sqrt{n}\left(\hat{\gamma}_n - \gamma_0 \right)$ is asymptotically normal.

(a) Consider the term $I_{n1}^{(k)} (k = 1,2)$.  If follows from Lemma \ref{lemmaA1} below that, for any fixed $0 < b < \infty$, there exists a positive constant $\epsilon>0$ such that, for $k =1, 2$,
\begin{flalign*}
\sup_{x \in R,\, \|\vartheta\| \le b}\sqrt{n}\big | H_{nk}(x;\vartheta) - H_{nk}(x;0) - H_k(x;\vartheta ) + H_k(x;0)\big| = O_p(n^{-\epsilon}),
\end{flalign*}
and therefore
\begin{flalign*}
\sup_{x \in R}|H_{nk}(x;\hat{\vartheta}_n) - H_{nk}(x;0) - H_k(x;\hat{\vartheta}_n) + H_k(x;0)| = O_p\big(n^{-1/2 - \epsilon}\big).
\end{flalign*}
Note that, through integration by parts, for a function $G(x)$,
\begin{flalign*}
\int_{n^{-2}}^{n^2}\log|x| dG(x) = \log(x)G(x)|_{n^{-2}}^{n^2} - \int_{n^{-2}}^{n^2} G(x) d\log(x).
\end{flalign*}
Hence,
\begin{flalign*}
\Big|\int_{n^{-2}}^{n^2}\log|x| dG(x)\Big| \le 8 \log(n) \sup_{x \in [n^{-2},n^2]} |G(x)|.
\end{flalign*}
Similarly,
$\big|\int_{-n^{2}}^{-n^{-2}}\log|x| dG(x)\big| \le 8 (\log n) \sup_{x \in [-n^{2}, \,-n^{-2}]} |G(x)|.$
Applying them to  $ |I_{n1}^{(k)}| $, we immediately have that, for $k =1,2$,
\begin{flalign*}
 |I_{n1}^{(k)}| &\le 16 \log(n) \sup_{x \in R}\big | H_{nk}(x;\hat{\vartheta}_n ) - H_{nk}(x;0) - H_k(x;\hat{\vartheta}_n ) + H_k(x;0)\big|\\
 & = o_p(n^{-1/2}).
\end{flalign*}

(b) Consider $I_{n2}^{(k)} (k = 1,2)$. Observe that
$\sqrt{n} \big (1- P (|\eta_t \sqrt{\alpha_0} + \phi_0| \in \mathcal{I}_n) \big) = o(1)$.
Thus,
\begin{flalign*}
\sqrt{n} I_{n1}^{(k)} &=  \frac{1}{n} \sum_{t=1}^n \big (I(|\eta_t \sqrt{\alpha_0}  - (-1)^k\phi_0| \in \mathcal{I}_n) \log |\eta_t \sqrt{\alpha_0} - (-1)^k\phi_0 |- \int_{\mathcal{I}_n} \log |x| dH_k(x;0)\big)\\
&= \frac{1}{\sqrt{n}} \sum_{t=1}^n \zeta_{kt} + o_p(1),
\end{flalign*}
and consequently,
\begin{eqnarray*}
\sqrt{n}\big(I_{n1}^{(1)} + I_{n2}^{(2)}\big)=\frac{1}{\sqrt{n}}\sum_{t=1}^n(\zeta_{1t} + \zeta_{2t}) + o_p(1).
\end{eqnarray*}

(c) Consider the term $I_{n3}^{(k)}$ $(k=1,2)$. Since $E \eta^2 < \infty$ and $f(\cdot)$ is bounded,
$$
|I_{n3}^{(k)}|^2 \le {E \big(\log |(-1)^{k-1}\phi_0 + \eta \sqrt{\alpha_0}|\big)^2 P (|(-1)^{k-1}\phi_0 + \eta \sqrt{\alpha_0}| \notin \mathcal{I}_n) }= o(n^{-1}).
$$

(d) Consider the term $I_{n4}^{(k)}$ $(k=1,2)$.
Define
$$
\tilde{v}_{t}^{(k)} = (-1)^{k-1} \Big( 1 + \frac{ (-1)^{k}y_{t-1}\sqrt{\alpha_0}}{\sqrt{{\omega}_0 + {\alpha}_0 y_{t-1}^2}} \Big),\quad
\tilde{u}_{1t} =\frac{1}{\sqrt{\alpha_0}}\tilde{u}_{2t},\quad
\tilde{u}_{2t} =\frac{-\sqrt{\alpha_0}}{2(\omega_0 + \alpha_0 y_{t-1}^2)}.
$$
By Taylor's expansion, $v_{t}^{(k)}(\hat{\theta}_n) (k =1,2)$ and $u_{t}(\hat{\theta}_n) $ can be approximated by
\begin{flalign*}
v_{t}^{(k)}(\hat{\theta}_n) &= \tilde{v}_{t}^{(k)} (\hat{\phi}_n - \phi_0) + o_p(n^{-1/2}),\\
u_{t}(\hat{\theta}_n) &= \tilde{u}_{1t} ({\hat{\alpha}_n} - {\alpha_0}) + \tilde{u}_{2t} ( \hat{\omega}_n - \omega_0 ) +o_p(n^{-1/2}),
\end{flalign*}
respectively. It follows from the laws of large numbers that
\begin{flalign*}
\frac{1}{n}\sum_{t = 1}^n \tilde{u}_{kt} \to \nu_k, a.s., \quad k = 1,2.
\end{flalign*}
Hence, it follows from Lemma \ref{lemmaA2} and Theorem \ref{theorem2} that
\begin{flalign*}
\sqrt{n} \big( I_{n4}^{(1)} + I_{n4}^{(2)} \big ) &= 2 \mu_1 \sqrt{n}(\hat{\phi}_n - \phi_0 )  + 2 \mu_2 \nu_1 \sqrt{n}({\hat{\alpha}_n} - {\alpha_0})+ 2\mu_2 \nu_2 \sqrt{n}( \hat{\omega}_n - \omega_0 ) +o_p(1)\\
& = \frac{1}{\sqrt{n}} \sum_{t=1}^n \big ( 2 \mu_1 \xi_{1t} + 2 \mu_2 \nu^T \xi_{2t} \big ) + o_p(1).
\end{flalign*}
The part (d) follows.

(II) Consider the nonstationary case, i.e., $\gamma_0 > 0$. Let $v_{nt}^{(k)}(\vartheta ) = n^{1/2} v_{t}^{(k)}\big((\phi_0 + n^{-1/2} \vartheta_1, \alpha_0 + n^{-1/2} \vartheta_2, \omega_0 + \vartheta_3)^T\big)$ and $u_{nt}(\vartheta ) = n^{1/2} u_{t}\big((\phi_0 + n^{-1/2} \vartheta_1, \alpha_0 + n^{-1/2} \vartheta_2, \omega_0 + \vartheta_3)^T \big)$. Define
\begin{flalign*}
\tilde{H}_{nk}(x;\vartheta) &= \frac{1}{n}\sum_{t = 1}^n I\Big( \eta_t  \le \frac{x + (-1)^{k}\phi_0 - n^{-1/2}v_{nt}^{(k)}(\vartheta)}{\sqrt{\alpha_0} + n^{-1/2}u_{nt}(\vartheta)}\Big),\\
\tilde{H}_k(x;\vartheta) &= \frac{1}{n}\sum_{t = 1}^n F\Big(\frac{x + (-1)^{k}\phi_0 - n^{-1/2}v_{nt}^{(k)}(\vartheta)}{\sqrt{\alpha_0} + n^{-1/2} u_{nt}(\vartheta)}\Big).
\end{flalign*}
Denote $\hat{\vartheta}_n = \big( n^{1/2}(\hat{\phi}_n - \phi_0), n^{1/2}(\hat{\alpha}_n - \alpha_0), \hat{\omega}_n - \omega_0\big)^T$,  then
\begin{flalign*}
2\hat{\gamma}_n =  \int_{\mathcal{I}_n} \log |x| d \tilde{H}_{n1}(x;\hat{\vartheta}_n) + \int_{\mathcal{I}_n} \log |x| d \tilde{H}_{n2}(x;\hat{\vartheta}_n) = \hat{\gamma}_{n1} + \hat{\gamma}_{n2}.
\end{flalign*}
We decompose the $\hat{\gamma}_{nk} - \gamma_0$ $(k = 1,2)$, respectively, as
\begin{flalign*}
\hat{\gamma}_{nk} - \gamma_0 &=
\int_{\mathcal{I}_n} \log|x| d \left( \tilde{H}_{nk}(x;\hat{\vartheta}_n) - \tilde{H}_{nk}(x;0) - \tilde{H}_k(x;\hat{\vartheta}_n ) + \tilde{H}_k(x;0)\right)\\
&\quad + \int_{\mathcal{I}_n} \log |x| d(\tilde{H}_{nk}(x;0) - \tilde{H}_k(x;0)) + \int_{R/\mathcal{I}_n} \log |x| d \tilde{H}_k(x;0)  \\
&\quad + \int_{\mathcal{I}_n} \log |x| d(\tilde{H}_k(x;\hat{\vartheta}_n) - \tilde{H}_k(x;0))\\
&:= I_{n1}^{(k)} + I_{n2}^{(k)} + I_{n3}^{(k)} + I_{n4}^{(k)}.
\end{flalign*}
To complete the proof, it suffices to show that
\begin{flalign*}
&(a)\quad \sqrt{n} I_{n1}^{(k)} = o_p(1),\quad k = 1,2;\\
&(b)\quad
\sqrt{n} \big ( I_{n2}^{(1)} + I_{n2}^{(2)} \big )= \frac{1}{\sqrt{n}} \sum_{t=1}^n
\big (\zeta_{1t} + \zeta_{2t}\big ) + o_p(1);\\
&(c)\quad \sqrt{n} I_{n3}^{(k)} = o_p(1),\quad k = 1,2; \\
&(d)\quad
\sqrt{n} \big ( I_{n4}^{(1)} + I_{n3}^{(2)} \big ) = \frac{1}{\sqrt{n}} \sum_{t=1}^n
\big ( 2\mu_1 \xi_{1t}\big)+ o_p(1).&&
\end{flalign*}
Combining (a)-(d), we immediately obtain that
$$
\sqrt{n}\left(\hat{\gamma}_n - \gamma_0 \right) = \frac{1}{\sqrt{n}} \sum_{t=1}^n \big(
\zeta_{t} + \mu_1 \tilde{\xi}_{1t}\big) + o_p(1).
$$
In fact, (a), (b) and (c) can be obtained in a similar fashion as those in the stationary case.  Now we only need to prove (d).
First, it follows from Lemma \ref{lemmaA1} below and Taylor's expansion that
\begin{flalign*}
\sqrt{n} \big( I_{n4}^{(1)} + I_{n4}^{(2)} \big ) &=   \frac{\mu_1}{\sqrt{n}}\sum_{t=1}^n \big \{ v_{t}^{(1)}(\hat{\theta}_n)+ v_t^{(2)}(\hat{\theta}_n) \big \}+  \frac{2\mu_2}{\sqrt{n}}\sum_{t=1}^n  u_{t}(\hat{\theta}_n) +o_p(1)\\
& = \frac{1}{\sqrt{n}} \sum_{t=1}^n \big ( 2 \mu_1 \tilde{\xi}_{1t} \big )+
\frac{2\mu_2}{\sqrt{n}}\sum_{t=1}^n  u_{t}(\hat{\theta}_n) + o_p(1).
\end{flalign*}
By the compactness of $\Theta$ and the fact (\ref{fact}), we have
\begin{eqnarray*}
\sum_{t=1}^n |u_{t}(\hat{\theta}_n)|\leq M\sum_{t=1}^n \frac{1}{1 + y_{t-1}^2}<\infty.
\end{eqnarray*}
Thus, $\sum_{t=1}^n |u_{t}(\hat{\theta}_n)|=O_p(1)$ and
\begin{eqnarray*}
\sqrt{n} \big( I_{n4}^{(1)} + I_{n4}^{(2)} \big )  = \frac{1}{\sqrt{n}} \sum_{t=1}^n \big ( 2 \mu_1 \tilde{\xi}_{1t} \big ) + o_p(1).
\end{eqnarray*}
The proof is complete.

\subsection{Proof of Theorem \ref{theorem5}}
Consider the case of $\gamma_0 < 0$. Using similar arguments in the proof of Theorem \ref{theorem4}, we obtain that
\begin{eqnarray*}
\sqrt{n} \left ( \hat{\gamma}^*_{n} - \gamma_0 \right ) = \frac{1}{\sqrt{n}}\sum_{t = 1}^n \varpi_t\big (  \zeta_t  +  \mu_1 \xi_{1t} + \mu_2 \nu^T \xi_{2t} \big ) + o_{P^*}(1).
\end{eqnarray*}
Then, with the asymptotic representation of $\sqrt{n}(\hat{\gamma}_n - \gamma_0)$, it follows that
\begin{eqnarray*}
\sqrt{n} \left ( \hat{\gamma}^*_{n} - \hat{\gamma}_n \right ) = \frac{1}{\sqrt{n}}\sum_{t = 1}^n
(\varpi_t-1) \big ( \zeta_{t} + \mu_1 \xi_{1t} + \mu_2 \nu^T \xi_{2t} \big ) + o_{P^*}(1).
\end{eqnarray*}
Thus, conditionally on $\{y_0,y_1,\cdots,y_n\}$, in probability, as $n \to \infty$,
\begin{eqnarray*}
\sqrt{n}\big (\hat{\gamma}_n^* - \hat{\gamma}_n \big )  \stackrel{\mathcal{L}}{\longrightarrow} N\big(0, ~ E (\zeta_t + \mu_1 \xi_{1t} + \mu_2 \nu^T \xi_{2t})^2\big).
\end{eqnarray*}
The results for $\gamma_0 > 0$ can be proved analogously. The proof is complete.

\subsection{Proof of Theorem \ref{theorem6}} The proof is immediate from Theorem \ref{theorem5}, and hence it is omitted.

\section{Lemmas with proofs}
In this subsection, two lemmas are given for the proofs of Theorems \ref{theorem4} and \ref{theorem5}.
\begin{lemma}\label{lemmaA1}
Let $\{\eta_t\}$ be i.i.d. with the distribution $F(\cdot)$, which has an a.e. positive bounded density $f(\cdot)$ with $\int_{x \in R} |x|f(x) dx < \infty$. Let $\mathcal{A}_{ni}$ be an array of sub $\sigma$-fields such that $\mathcal{A}_{nt} \subset \mathcal{A}_{n(t+1)}, 1 \le t \le n,n \ge 1$; for each $\theta$, $(v_{n1}(\theta),u_{n1}(\theta))$ is $\mathcal{A}_{n1}$-measurable and $(\eta_1,\cdots,\eta_{t-1},v_{nt}(\theta),u_{nt}(\theta))$, $t \le j$ are $\mathcal{A}_{nj}$ measurable, $2 \le j  \le n$; $\eta_t$ is independent of $\mathcal{A}_{nt}, 1 \le t \le n$. For $v_{nt}(\theta)$ and $u_{nt}(\theta)$, assume that
\begin{eqnarray*}
\frac{1}{n}\sum_{t = 1}^n \sup_{||\theta_1 - \theta_2|| \le \delta} \Big|v_{nt}(\theta_1) - v_{nt}(\theta_2) \Big| \le C_1 \delta, \quad a.s.,
\\
\mbox{and} \quad
\frac{1}{n}\sum_{t = 1}^n \sup_{||\theta_1 - \theta_2|| \le \delta} \Big|u_{nt}(\theta_1) - u_{nt}(\theta_2) \Big| \le C_2 \delta, \quad a.s.
\end{eqnarray*}
for some universal positive constants $C_1$ and $C_2$. For each $x \in R$, define
\begin{flalign*}
H_n(x;\theta) =\frac{1}{n}\sum_{t = 1}^n I\Big( \eta_t  \le \frac{x  - n^{-1/2}v_{nt}(\theta)}{\sqrt{\alpha_0} + n^{-1/2}u_{nt}(\theta)}\Big)\quad\mbox{and} \quad
H(x;\theta) =\frac{1}{n}\sum_{t = 1}^n F\Big(\frac{x - n^{-1/2}v_{nt}(\theta)}{\sqrt{\alpha_0} + n^{-1/2} u_{nt}(\theta)}\Big),
\end{flalign*}
where $\alpha_0$ are a positive constant.
Then, for any constant $b$ with $0 < b < \infty$, there exists a positive constant $\epsilon$ such that
\begin{flalign*}
\sup_{x \in R, \:\|\theta\| \le b}\sqrt{n}\left | H_n(x;\theta) - H_n(x;0) - H(x;\theta) + H(x;0)\right| = O_p(n^{-\epsilon}).
\end{flalign*}
\end{lemma}
{\sc Proof.} This lemma is similar to that of Lemma 8.3.2 in \cite{koul2002} but the rate of convergence is strengthened. It should be noted that we use the pseudometric
\begin{flalign*}
d_b(x,y) = \sup_{|z| \le b} \Big|F\Big(\frac{x - \phi_0 - n^{-1/2} z_1}{ \sqrt{\alpha_0} + n^{-1/2} z_2}\Big) - F\Big(\frac{y - \phi_0 - n^{-1/2} z_1}{ \sqrt{\alpha_0} + n^{-1/2} z_2}\Big) \Big|
\end{flalign*}
for this local and scale setup, where $z = (z_1,z_2) \in R^2$, $|z| = |z_1| \vee |z_2|$, $x,y \in R$, $b > 0$. If we let $\mathfrak{N}(\delta,b)$ be the cardinality of the minimal $\delta$-net of $(R,d_b)$, we can show that
$ \mathfrak{N}(\delta,b) \le C_b \delta^{-4} $ with  $C_b < \infty$ for any $0<n^{-1/2}b<1$. We can repeat the similar arguments as in Lemma 8.3.2 in \cite{koul2002} using the metric $d_b$. The convergence rate is strengthened since we make stronger conditions of the terms $u_{nt}(\theta)$ and $v_{nt}(\theta)$. The details are omitted here.

\begin{lemma}\label{lemmaA2}
Suppose that $\eta$ is a random variable with the density $f(\cdot)$ and the derivative $f'(x)$ of $f(x)$ exists for each $x$. The constants $\phi_0$ and $\alpha_0$ are finite and $\alpha_0 >0$.
Assume $b_1 = - \alpha_0^{-1/2}\int \log|\phi_0 + \eta \sqrt{\alpha_0}|f'(\eta)d \eta < \infty$ and
$b_2 = - \alpha_0^{-1/2}\int \log|\phi_0 + \eta \sqrt{\alpha_0}| (f(\eta) + \eta f'(\eta )  ) d \eta < \infty.$
Denote $\delta = (\delta_1,\delta_2)^T$ and
$ h(\eta;\delta) = \log |\phi_0 + \eta \sqrt{\alpha_0} + \delta_1 + \eta \delta_2|.$
Then, for small constants $\delta_1$ and $\delta_2$ with $\delta_1 \to 0$ and $\delta_2 \to 0$, the following approximation holds:
\begin{eqnarray*}
E \big(h(\eta;\delta)-h(\eta;0)\big)=b_1\delta_1+b_2\delta_2+o(|\delta_1|+|\delta_2|).
\end{eqnarray*}
\end{lemma}

{\sc Proof.} $E h(\eta;\delta)$ is expressed as
$$
Eh(\eta;\delta) = \int  f\big(\frac{\eta - \phi_0 - \delta_1}{\sqrt{\alpha_0}+ \delta_2}\big) \frac{\log|\eta|}{{\sqrt{\alpha_0}+ \delta_2}}d \eta.
$$
Hence,
\begin{flalign*}
E \big(h(\eta;\delta) - h(\eta;0)\big)
&= \int \log|\eta| \Big\{ f\big(\frac{\eta - \phi_0 - \delta_1}{\sqrt{\alpha_0} + \delta_2} \big) - f\big(\frac{\eta - \phi_0}{\sqrt{\alpha_0}} \big) \Big\} d \eta \frac{1}{{\sqrt{\alpha_0}+ \delta_2}}\\
& \hskip 1cm + \Big( \frac{1}{{\sqrt{\alpha_0}+ \delta_2}} - \frac{1}{{\sqrt{\alpha_0}}} \Big) \int \log|\eta| f\big(\frac{\eta - \phi_0}{\sqrt{\alpha_0}} \big) d \eta
\\
& = I_1 + I_2.
\end{flalign*}
Let us handle the first term $I_1$.
By Taylor 's expansion for $f(\cdot)$ at the point $(\eta - \phi_0)/\sqrt{\alpha_0}$, we obtain that
\begin{flalign*}
f\big(\frac{\eta - \phi_0 - \delta_1}{\sqrt{\alpha_0}+ \delta_2} \big) - f\big(\frac{\eta - \phi_0}{\sqrt{\alpha_0}} \big) 
&= f'\big(\frac{\eta - \phi_0}{\sqrt{\alpha_0}} \big) \Big(\frac{\eta - \phi_0 - \delta_1}{\sqrt{\alpha_0}+ \delta_2} - \frac{\eta - \phi_0}{\sqrt{\alpha_0}}\Big) +o(|\delta_1| + |\delta_2|) \\
&= f'\big(\frac{\eta - \phi_0}{\sqrt{\alpha_0}} \big) \Big(\frac{ - \delta_1 \sqrt{\alpha_0} - \delta_2\eta + \phi_0 \delta_2}{\alpha_0} \Big) +o(|\delta_1| + |\delta_2|).
\end{flalign*}
Thus,
\begin{flalign*}
I_1 &=  \int \log|\eta|f'\big(\frac{\eta - \phi_0}{\sqrt{\alpha_0}} \big) \Big(\frac{ - \delta_1 \sqrt{\alpha_0} - \delta_2\eta + \phi_0 \delta_2}{\alpha_0} \Big)d \eta  \frac{1 }{\sqrt{\alpha_0}}+ o(|\delta_1| + |\delta_2|) \\
&=- \alpha_0^{-1/2}\delta_1\int \log|\phi_0 + \eta \sqrt{\alpha_0}|f'(\eta) d \eta   \\
&\quad-\alpha_0^{-1/2}\delta_2 \int \log|\phi_0 + \eta \sqrt{\alpha_0}| f'(\eta) \eta d \eta +o(|\delta_1| + |\delta_2|).
\end{flalign*}
Similarly, we have that
\begin{flalign*}
I_2=- \alpha_0^{-1/2}\delta_2 \int \log|\phi_0 + \eta \sqrt{\alpha_0}| f\left(\eta \right ) d \eta + o(\delta_2).
\end{flalign*}
Combining the above results, the proof follows.

\end{document}